\documentclass[hidelinks,onefignum,onetabnum]{siamart250211}

\usepackage{lipsum}
\usepackage{bm}
\usepackage{amsfonts}
\usepackage{graphicx}
\usepackage{epstopdf}
\usepackage{subfigure}
\usepackage{algorithmic}
\ifpdf
  \DeclareGraphicsExtensions{.eps,.pdf,.png,.jpg}
\else
  \DeclareGraphicsExtensions{.eps}
\fi

\newsiamremark{remark}{Remark}
\newsiamremark{hypothesis}{Hypothesis}
\crefname{hypothesis}{Hypothesis}{Hypotheses}
\newsiamthm{claim}{Claim}
\newsiamremark{fact}{Fact}
\crefname{fact}{Fact}{Facts}

\headers{ROM of iterative decoupled algorithms for Biot's model}{H. Gu, F. Ballarin, M. Cai, and J. Li}

\title{
POD-based reduced order modeling of global-in-time iterative decoupled algorithms for Biot's consolidation model\thanks{Submitted to arXiv.
\funding{H. Gu is supported by National NSF of China No. 123B2016. F. Ballarin acknowledges the PRIN 2022 PNRR project ``ROMEU: Reduced Order Models for Environmental and Urban flows'' funded by the European Union -- NextGenerationEU under the National Recovery and Resilience Plan (NRRP), Mission 4 Component 2, CUP J53D23015960001. M. Cai is supported in part by the NIH-RCMI grant through U54MD013376, and the affiliated project award from the Center for Equitable Artificial Intelligence and Machine Learning Systems (CEAMLS) at Morgan State University (project ID 02232301). J. Li is supported by the Shenzhen Sci-Tech FundRCJC20200714114556020, Guangdong Basic and Applied Research Fund 2023B1515250005.
}
\\
\and \textcolor{white}{HP} \thanks{A} \ School of Humanities and Fundamental Sciences, Shenzhen Institute of Information Technology, Shenzhen, 518172, China. Email: \email{guhuipeng@sziit.edu.cn}
\\
\and \textcolor{white}{FB} \thanks{A} \ Department of Mathematics and Physics, Università Cattolica del Sacro Cuore, Brescia, 25133, Italy. Email:
\email{francesco.ballarin@unicatt.it}.
\\
\and \textcolor{white}{fB} \thanks{A} \ Department of Mathematics for Economic, Financial and Actuarial Sciences, Università Cattolica del Sacro Cuore, Milano, 20123, Italy.
\\
\and \textcolor{white}{C} \thanks{A} \ Department of Mathematics, Morgan State University, Baltimore, MD 21251, USA. Email:
\email{mingchao.cai@morgan.edu}
\\
\and \textcolor{white}{JZ} \thanks{A} \ Department of Mathematics, Southern University of Science and Technology, Shenzhen, 518055, China. Email: \email{li.jz@sustech.edu.cn}
}
}

\author{
Huipeng Gu\footnotemark[2] ,
Francesco Ballarin\footnotemark[3] $^,$\footnotemark[4]  ,
Mingchao Cai\footnotemark[5] ,
Jingzhi Li\footnotemark[6]
}

\usepackage{amsopn}

\ifpdf
\hypersetup{
  pdftitle={POD-based reduced order modeling of global-in-time iterative decoupled algorithms for Biot's consolidation model},
  pdfauthor={H. Gu, F. Ballarin, M. Cai, and J. Li}
}
\fi

\begin{document}

\maketitle

\begin{abstract}
This paper focuses on the efficient numerical algorithms of a three-field Biot's consolidation model. 
The approach begins with the introduction of innovative monolithic and global-in-time iterative decoupled algorithms, which incorporate the backward differentiation formulas for time discretization.
In each iteration, these algorithms involve solving a diffusion subproblem over the entire temporal domain, followed by solving a generalized Stokes subproblem over the same time interval.
To accelerate the global-in-time iterative process, we present a reduced order modeling approach based on proper orthogonal decomposition, aimed at reducing the primary computational cost from the generalized Stokes subproblem. 
The effectiveness of this novel method is validated both theoretically and through numerical experiments.
\end{abstract}



\section{Introduction}
Biot's consolidation model, which describes the interaction between fluid flow and mechanical deformation in porous media, has broad applications in fields such as geomechanics, petroleum engineering, and biomechanics \cite{biot1941general, ju2020parameter}. 
Due to the complexity of real-world conditions, obtaining analytical solutions for this coupled system is highly challenging.
Consequently, researchers have developed numerical methods to approximate solutions, including approaches based on finite element methods (FEMs), finite volume methods (FVMs), and other discretization techniques \cite{phillips2007coupling,naumovich2006finite,cai2023combination}. 

To address this time-dependent multi-physics system, monolithic numerical methods are employed to solve the entire coupled system directly, providing both accuracy and stability.
However, these methods are computationally expensive, particularly when applied to problems with extensive temporal domains or high spatial resolutions.
In order to overcome this challenge, a common strategy for reducing computational costs is to decouple the system into smaller subproblems (separate the mechanics equation
and the flow equation). 
This approach can be broadly classified into two categories: non-iterative \cite{chaabane2018splitting, feng2018analysis, altmann2021semi, cai2023some} and iterative decoupled algorithms \cite{kim2011stability, mikelic2013convergence, mikelic2014numerical,gu2023iterative}.
Non-iterative decoupled schemes typically employ time-delay techniques that significantly reduce computational costs, but they may compromise stability, particularly when higher-order methods (in time) are used.
Iterative decoupled schemes, such as fixed-stress methods, provide a balanced alternative. They enhance stability in comparison to non-iterative approaches while simultaneously reducing computational costs compared to monolithic methods, thus making them a favorable option for achieving both efficiency and accuracy.

Among iterative decoupled approaches, several innovative strategies have been proposed to further reduce computational costs and time. Multirate iterative schemes \cite{almani2016convergence, almani2023convergence} achieve significant efficiency by employing multiple finer time steps for the flow problem within a single coarse time step for the mechanics problem. Multiscale iterative schemes \cite{dana2018multiscale} adopt a similar efficiency-driven strategy, solving the flow problem on a finer mesh, while solving the mechanics problem on a coarser mesh. Global-in-time iterative schemes \cite{borregales2019partially, ahmed2020adaptive, gu2024crank} take a different approach, leveraging parallel computing to solve the mechanics problem across the entire temporal domain, thus enabling efficient handling of time-dependent processes.
In addition to these strategies, recent advancements in reduced order modeling (ROM) techniques have shown considerable potential for further reducing computational burdens. 
For instance, proper orthogonal decomposition (POD) is employed to construct a reduced basis that is then employed in the definition of a reduced order model via Galerkin projection
\cite{ballarin2024projection}, enhancing computational efficiency.
Furthermore, the dual-weighted residual method has been utilized to develop error-controlled incremental POD-based ROMs \cite{fischer2024more}, making them well-suited for real-time simulations.

In this work, we propose a POD-based reduced order model for global-in-time iterative decoupled algorithms using the backward differentiation formula. Applied to a three-field formulation of Biot's consolidation model \cite{oyarzua2016locking}, the flow equation is treated as a diffusion problem, while the mechanics equation is formulated as a generalized Stokes problem. 
When solving the generalized Stokes problem in each iteration, the approach involves using the computationally expensive full order model (FOM) for only a small portion of the entire time interval. These solutions serve as snapshots to generate the reduced basis (RB) using the POD method, enabling the reconstruction of numerical solutions for the entire time interval at a significantly reduced cost. Since the principal challenge arises from the resolution of the mechanical equation, this approach effectively reduces the computational burden while maintaining accuracy. 
Additionally, we consider backward differentiation formulas (BDF) for time discretization \cite{fu2019high, altmann2024higher}, focusing on the second-order BDF scheme due to its unconditional stability.
Moreover, we demonstrate the efficiency of the proposed method both theoretically and through numerical experiments, highlighting its ability to reduce computational costs while maintaining solution accuracy. 

The paper is structured as follows: \Cref{Sec2} presents the model formulation, while \Cref{Sec3} discusses the fully discrete monolithic algorithm and the global-in-time iterative decoupled algorithm based on BDF, along with its convergence analysis. \Cref{Sec4} introduces the POD-based reduced order modeling of global-in-time iterative schemes and provides the theoretical analysis. Finally, \Cref{Sec5} presents numerical tests to demonstrate the performance and efficiency of the proposed approach.

\section{Model formulation}
\label{Sec2}
Let $\Omega \subset \mathbb{R}^d$ ($d= 2$ or $3$) be a bounded polygonal domain with boundary $\partial \Omega$, and let $T < \infty$ denote the final time. We denote by $H^1(\Omega)$ the classical Sobolev space of functions whose distributional derivatives up to first order belong to $L^2(\Omega)$, i.e., 
\begin{align*}
    H^1(\Omega) = \{v \in L^2(\Omega); \nabla v \in (L^2(\Omega))^d \},
\end{align*}
equipped with the norm $\|\cdot\|_{H^1}$ and semi-norm $|\cdot|_{H^1}$. For simplicity, $(\cdot, \cdot)$ and $(\cdot, \cdot)_1$ denote the $L^2$ inner product and $H^1$ inner product, respectively.

In this study, we focus on the three-field formulation of Biot's consolidation model, as presented in the references \cite{oyarzua2016locking,lee2017parameter}. This formulation introduces an intermediate variable $\xi$, called the total pressure, as one of the unknowns. The system is described by the following equations:
\begin{align}
    - \nabla \cdot ( 2\mu \varepsilon (\bm{u}) - \xi \mathbf{I} ) & = \bm{f} \quad \ \  \mbox{in} \ \Omega \times (0,T], \label{threefield1} 
    \\
    \nabla \cdot \bm{u} + \frac{1}{\lambda} \xi - \frac{\alpha}{\lambda} p & = 0 \quad \quad \mbox{in} \ \Omega \times (0,T], \label{threefield2} 
    \\
    \left( c_0 + \frac{\alpha^2}{\lambda} \right) \partial_t p - \frac{\alpha}{\lambda} \partial_t \xi - \nabla \cdot( k_p \nabla p ) & = g \quad \quad \mbox{in} \ \Omega \times (0,T].  \label{threefield3} 
\end{align}
Here, the primary unknowns are the displacement vector of the solid $\bm{u}$, the total pressure $\xi$ and the fluid pressure $p$. The symbol $\varepsilon(\bm{u}) = \frac{1}{2} [ \nabla \bm{u} + (\nabla \bm{u})^T ]$ is the linearized strain tensor, $\alpha > 0$ is the Biot-Willis constant which is close to 1, $\bm{f}$ is the body force, $c_0 \geq 0$ is the specific storage coefficient, $k_p$ represents the hydraulic conductivity, $g$ is a source or sink term, $\mathbf{I}$ is the identity matrix, and Lam{\'e} constants $\lambda$ and $\mu$ are computed from the Young's modulus $E$ and the Poisson ratio $\nu$:
\begin{align*}
	\lambda = \frac{E\nu}{(1+\nu)(1-2\nu)},\ \ \ \mu = \frac{E}{2(1+\nu)}.
\end{align*}
Appropriate initial and boundary conditions must be specified to complete the system of governing equations. We assume that the boundary $\partial \Omega$ can be divided 
as follows:
\begin{align}
    \partial \Omega = \Gamma_u \cup \Gamma_\sigma = \Gamma_p \cup \Gamma_q, \ \text{where} \ |\Gamma_u|, \ |\Gamma_\sigma|, \ |\Gamma_p|, \ \text{and} \ |\Gamma_q|>0, 
    \label{Boundarysets}
\end{align}
with $\Gamma_u \cap \Gamma_\sigma = \emptyset$ and $\Gamma_p \cap \Gamma_q = \emptyset$.
We consider the mixed-type boundary conditions:
\begin{align}
    \bm{u} = \bm{0} \quad & \mbox{on} \ \Gamma_{u},  \label{bc1} \\
    \left( 2 \mu \varepsilon(\bm{u}) - \xi \mathbf{I} \right) \bm{n} = \bm{0} \quad & \mbox{on} \ \Gamma_{\sigma}, \\ 
    p=0 \quad & \mbox{on} \ \Gamma_p, \\ 
    k_p (\nabla p ) \cdot \bm{n} = 0 \quad & \mbox{on} \ \Gamma_q,  \label{bc4}
\end{align}
along with the following initial conditions:
\begin{align}
    \bm{u}(0) = \bm{u}^0, \quad p(0) = p^0, \quad \xi(0) = \alpha p^0 - \lambda \nabla \cdot \bm{u}^0. \label{ic}
\end{align}

Here, we define the function spaces
\begin{align*}
    \bm{V} = \bm{H}^1_{0,\Gamma_{u}} (\Omega) = \left[ H^1_{0,\Gamma_{u}} (\Omega) \right]^d, \quad W = L^2(\Omega), \quad M = H^1_{0,\Gamma_p} (\Omega),
\end{align*}
where $H^1_{0,\Gamma} = \{ u \in H^1(\Omega) : u|_{\Gamma} = 0 \}$,
and the corresponding bilinear forms as
    \begin{align*}
        & a_1(\bm{u},\bm{v}) = 2\mu \int_\Omega  \varepsilon (\bm{u}) :  \varepsilon (\bm{v}),  \ \ \ b(\bm{v},w) =  \int_\Omega w \  \nabla \cdot \bm{v}, \\
        & a_2(\xi,w) = \frac{1}{\lambda} \int_\Omega \xi w,  \quad \quad \quad \quad \ \ c(p , w) =  \frac{\alpha}{\lambda} \int_\Omega p w, \\
        & a_3(p,q) = \Big( c_0+\frac{\alpha^2}{\lambda} \Big) \int_\Omega p q, \ \ \ \ d(p,q) =  k_p \int_\Omega \nabla p \cdot \nabla q,
    \end{align*}
where $\bm{u}, \bm{v} \in \bm{V}$, $\xi, w \in W$, $p, q \in M$.
This leads to the following variational formulation: For almost every $t \in (0,T]$, find $(\bm{u}, \xi, p) \in \bm{V} \times W \times M$ such that
    \begin{align}
        a_1(\bm{u},\bm{v})-b(\bm{v},\xi) & = ( \bm{f},\bm{v} ) ,  \label{WF1} 
        \\
        b(\bm{u},w) + a_2(\xi,w) - c(p,w) & = 0,  \label{WF2} 
        \\
        a_3(\partial_t p,q)-c(q,\partial_t \xi)+d(p,q) & = ( g, q ),  \label{WF3}
    \end{align}
for all $(\bm{v}, w, q) \in \bm{V} \times W \times M$. 
Furthermore, the following inf-sup condition \cite{girault1979finite, brenner1993nonconforming} holds: there exists a constant $\beta > 0$ depending only on $\Omega$ such that
\begin{align}
    \sup_{\bm{0} \not = \bm{v}\in \bm{V}} \frac{  b(\bm{v},w) }{\| \bm{v} \|_{H^1}} \geq \beta \| w \|_{L^2}, \ \ \ \forall w \in W. \label{infsupcon}
\end{align}
The well-posedness of problem \eqref{WF1}-\eqref{WF3} can be found in \cite{oyarzua2016locking}, and we assume that the solutions $\bm{u}$ and $p$, together with $\bm{f}$ and $g$, satisfy the required regularity conditions.

\section{Monolithic and global-in-time iterative decoupled algorithms for the fully discrete problem}
\label{Sec3}

In this section, we discuss about two fully discrete finite element methods of full order model (FOM) with backward differentiation formulas (BDF) for solving the problem \eqref{WF1}-\eqref{WF3} numerically: one based on monolithic algorithms and the other on global-in-time iterative decoupled algorithms.

\subsection{Spatial and Temporal Discretization}

For the temporal discretization, we partition the time interval $[0,T]$ into $N$ uniform subintervals with a fixed time step size $\Delta t$, such that $0 = t^0 < t^1 < \cdots < t^N = T$, where $t^n = n \Delta t$ for $n = 0, 1,\cdots,N$. The implicit method employed is the $m$-step backward differentiation formula (BDF-$m$) \cite{hairer2010solving}, which enables higher-order accuracy in time. The $m$-step backward difference operator is denoted as
\begin{align}
    \mathbb{D}_t^m y^n = \frac{1}{\Delta t} \Big( \sum_{j=0}^m \eta_j y^{n-j} \Big),
\end{align}
where $y^n:=y(t^n)$, and $\eta_j$ are the BDF-$m$ coefficients listed in Table \ref{Table:bdf}. 

\begin{table}[h]
\renewcommand\arraystretch{1.1}
\centering
\caption{BDF-$m$ coefficients $\eta_{j}$, for $1\leq m\leq 6$.}
\label{Table:bdf}
\begin{tabular}{c|ccccccc}
$m$ & $\eta_0$ & $\eta_1$ & $\eta_2$ & $\eta_3$ & $\eta_4$ & $\eta_5$ \\ \hline
1 & $1$ & $-1$ & & & & & \\
2 & $\frac{3}{2}$ & $-2$ & $\frac{1}{2}$ & & & & \\
3 & $\frac{11}{6}$ & $-3$ & $\frac{3}{2}$ & $-\frac{1}{3}$ & & & \\
4 & $\frac{25}{12}$ & $-4$ & $3$ & $-\frac{4}{3}$ & $\frac{1}{4}$ & & \\
5 & $\frac{137}{60}$ & $-5$ & $5$ & $-\frac{10}{3}$ & $\frac{5}{4}$ & $-\frac{1}{5}$ & \\
6 & $\frac{147}{60}$ & $-6$ & $\frac{15}{2}$ & $-\frac{20}{3}$ & $\frac{15}{4}$ & $-\frac{6}{5}$ & $\frac{1}{6}$ \\
\end{tabular}
\end{table}

As mentioned in \cite{suli2003introduction, fu2019high}, the BDF-$m$ methods are zero-stable only for $1\leq m \leq 6$, and are A-stable for $m=1$ and $m=2$.
In this work, we focus on the BDF-$2$ method due to its second-order accuracy and unconditional A-stability. For BDF-$m$ methods with $3 \leq m \leq 5$, the analysis can be similarly addressed using the multiplier technique \cite{nevanlinna1981multiplier}, while for BDF-$6$ the energy technique \cite{akrivis2021energy} is required, which is beyond the scope of this paper.
First, we introduce the $G$-matrix associated with the classical BDF-$2$ method \cite{chen2013efficient}, which is defined as
\begin{align*}
    G =
    \begin{bmatrix}
        \frac{1}{2} & -1 \\
        -1 & \frac{5}{2}
    \end{bmatrix},
\end{align*}
along with the corresponding $G$-norm,
$
    \|[y_1;y_2]\|_G^2 = \left([y_1;y_2] , G[y_1;y_2] \right)
$, which is equivalent to the $(L^2)^2$-norm: there exists $C_l, C_u >0$ such that
\begin{align*}
    C_l \|[y_1;y_2]\|_G \leq \|[y_1;y_2]\|_{L^2} \leq C_u \|[y_1;y_2]\|_G.
\end{align*}
Moreover, the following well-known equality holds:
\begin{align}
    ( \mathbb{D}_t^2 y^n , y^n ) = & \frac{1}{2 \Delta t} \|[y^n; y^{n-1}] \|_G^2 - \frac{1}{2 \Delta t} \|[y^{n-1}; y^{n-2}] \|_G^2 
    \nonumber
    \\ 
    & \ + \frac{\|y^n - 2y^{n-1} + y^{n-2}\|_{L^2}^2}{4 \Delta t}. \label{bdfineq}
\end{align}

For the spatial discretization, we partition the domain $\Omega$ into triangles in $\mathbb{R}^2$ or tetrahedra in $\mathbb{R}^3$, denoted as $\mathcal{T}_h$, where $h$ represents the maximum diameter of the elements in the mesh. The Taylor-Hood element pair $(\bm{V}_{h}, W_h)$ is employed for the variables $(\bm{u}, \xi)$, while the Lagrange finite element space $M_h$ is used for the variable $p$. These spaces are defined as follows:
\begin{align*} 
	& \bm{V}_{h} := \{ \bm{v}_h \in \bm{H}^1_{0, \Gamma_u} (\Omega) \cap {\bm{C} }^0(\bar{\Omega}); \ \bm{v}_h |_E  \in {\bm P}_{k}(E), ~\forall E \in \mathcal{T}_h \},
	\\
	& W_h := \{ w_h \in L^2(\Omega) \cap C^0(\bar{\Omega}); \ w_h |_E  \in  P_{k-1}(E), ~\forall E \in \mathcal{T}_h \},
	\\
	& M_h := \{ q_h \in H^1_{0, \Gamma_p} (\Omega) \cap C^0(\bar{\Omega}); \ q_h |_E  \in P_{l}(E), ~\forall E \in \mathcal{T}_h  \},
\end{align*}
where $k \geq 2$, $l \geq 1$ are integers. The Taylor-Hood element pair is known to satisfy the discrete inf-sup condition \cite{girault1979finite}: there exists a positive constant $\Tilde{\beta}$ independent of $h$, such that
\begin{align}
    \sup_{\bm{0} \not = \bm{v}_h \in \bm{V}_h} \frac{  b(\bm{v}_h,w_h) }{\| \bm{v}_h \|_{H^1}} \geq \tilde{\beta} \| w_h \|_{L^2}, \ \ \ \forall w_h \in W_h. \label{dinfsupcon}
\end{align}

\subsection{ Fully discrete finite element methods of full order model (FOM)} 

The monolithic method solves the coupled system \eqref{WF1}-\eqref{WF3} directly at each time step. This approach ensures strong coupling between all variables $(\bm{u},\xi,p)$, yielding a robust and accurate solution at the expense of solving a larger, fully coupled system.
The second-order fully discrete monolithic algorithm using the BDF-$2$ temporal discretization is presented in \cref{algo:Coupled_BDF2}.
\begin{algorithm}[H]
\raggedright
  \caption{: A fully coupled BDF-$2$ algorithm}
  \label{algo:Coupled_BDF2}
  \textbf{Input:} initial data $\{\xi_{h}^{n}\}_{n=0}^{1} \subset W_h$, $\{p_{h}^{n}\}_{n=0}^{1} \subset M_h$.
  \\
  \textbf{Output:} solutions $\{(\bm{u}_{h}^{n},\xi_{h}^{n},p_{h}^{n})\}_{n=2}^N \subset \bm{V}_h \times W_h \times M_h$.
  \\
  \textbf{for} $n$ from $2$ to $N$
  \\
  \quad \quad find $(\bm{u}_h^{n},\xi_h^{n},p_h^{n}) \in \bm{V}_h \times W_h \times M_h$ such that
  \vspace{-0.3cm}
  \begin{align}
    a_1(\bm{u}_h^{n},\bm{v}_h)-b(\bm{v}_h,\xi_h^{n} ) & = ( \bm{f}^{n},\bm{v}_h ) , \quad  \forall \bm{v}_h \in \bm{V}_h, \label{Coupled_BDF2_1}
    \\
    b(\bm{u}_h^{n},w_h) + a_2(\xi_h^{n},w_h) - c(p_h^{n},w_h ) & = 0
    , \quad \quad  \quad \quad  \forall w_h \in W_h, \label{Coupled_BDF2_2} 
    \\
    a_3 ( \mathbb{D}_{t}^2 p_h^{n} , q_h ) - c( q_h , \mathbb{D}_{t}^2 \xi_h^{n} ) 
    + d ( p_h^{n},q_h ) & =  ( g^{n} , q_h ), \quad \ \forall q_h \in M_h. \label{Coupled_BDF2_3}
  \end{align}
  \vspace{-0.7cm}
  \\
  \textbf{end for}
\end{algorithm}

To reduce the computational cost of solving the fully coupled system, the global-in-time iterative decoupled algorithms are introduced. These methods allow for the equations to be solved in a sequentially decoupled manner while iterating to achieve convergence for the global-in-time solution. The global-in-time BDF-$2$ algorithm is described in \cref{algo:GTIDA_bdf_fom}.

\begin{algorithm}[H]
\raggedright
  \caption{: A global-in-time iterative decoupled BDF-$2$ algorithm}
  \label{algo:GTIDA_bdf_fom}
  \textbf{Input:} initial data $\{\xi_{h,f}^{n}\}_{n=0}^{1} \subset W_h$ and $\{p_{h,f}^{n}\}_{n=0}^{1} \subset M_h$, initial guesses 
    \\
  \quad \quad \quad \ \ $\{\xi_{h,f}^{n,0}\}_{n=2}^{N} \subset W_h$, and initial iteration number $i = 0$.
  \\
  \textbf{Output:} solutions $\{(\bm{u}_{h,f}^{n},\xi_{h,f}^{n},p_{h,f}^{n})\}_{n=2}^N \subset \bm{V}_h \times W_h \times M_h$.
  \\
  \textbf{while} not converged \textbf{do}
  \\
  \quad \quad set $i = i + 1$,  $p_{h,f}^{0,i} = p_{h,f}^0$, $p_{h,f}^{1,i} = p_{h,f}^1$, $\xi_{h,f}^{0,i-1} = \xi_{h,f}^{0}$, and $\xi_{h,f}^{1,i-1} = \xi_{h,f}^{1}$. 
  \\
  \quad \quad \textbf{Step a:}  find $\{ p_{h,f}^{n,i} \}_{n=2}^{N} \subset M_h$ such that
  \vspace{-0.3cm}
  \begin{align}
    & a_3 ( \mathbb{D}_{t}^2 p_{h,f}^{n,i} , q_h ) 
    + d ( p_{h,f}^{n,i},q_h )  =  ( g^{n} , q_h ) + c( q_h , \mathbb{D}_{t}^2 \xi_{h,f}^{n,i-1} ) , \quad  \forall q_h \in M_h. \label{gt3}
  \end{align}
  \vspace{-0.7cm}
  \\
  \quad \quad \textbf{Step b:} find $\{(\bm{u}_{h,f}^{n,i},\xi_{h,f}^{n,i})\}_{n=2}^N \subset \bm{V}_h \times W_h$, such that
  \vspace{-0.3cm}
  \begin{align}
    & a_1(\bm{u}_{h,f}^{n,i},\bm{v}_h)-b(\bm{v}_h,\xi_{h,f}^{n,i}) = ( \bm{f}^{n},\bm{v}_h ), \quad \quad \ \forall \bm{v}_h \in \bm{V}_h, \label{gt1}
    \\
    & b(\bm{u}_{h,f}^{n,i},w_h) + a_2(\xi_{h,f}^{n,i},w_h) = c(p_{h,f}^{n,i},w_h)  ,  \quad   \forall w_h \in W_h. \label{gt2}
  \end{align}  
  \vspace{-0.7cm}
  \\
  \textbf{end while}
  \\
  set $(\bm{u}_{h,f}^{n},\xi_{h,f}^{n},p_{h,f}^{n}) = (\bm{u}_{h,f}^{n,i},\xi_{h,f}^{n,i},p_{h,f}^{n,i})$ for $n = 2, \cdots, N$.
\end{algorithm}

We note that \cref{algo:GTIDA_bdf_fom} requires additional initial guesses $\{\xi_{h,f}^{n,0}\}_{n=2}^{N} \subset W_h$, compared to \cref{algo:Coupled_BDF2}. 
In the absence of prior knowledge about the initial guesses, the first time step data $\xi_{h,f}^{1}$ can be used to initialize $\xi_{h,f}^{n,0}$ as $\xi_{h,f}^{n,0} = \xi_{h,f}^{1}$ for $n = 2, \dots, N$. The algorithm is checkpoint-friendly, enabling it to resume from stored values of $\{\xi_{h,f}^{n,i}\}_{n=2}^N$ at the $i$-th iteration. Furthermore, following the approach in \cite{gu2024crank}, we provide the convergence analysis of \cref{algo:GTIDA_bdf_fom}, as stated in \cref{Thm1}.

\begin{theorem}
\label{Thm1}
Let $(\bm{u}_h^n,\xi_h^n,p_h^n)$ and $(\bm{u}_{h,f}^{n,i},\xi_{h,f}^{n,i},p_{h,f}^{n,i})$ be the solutions of \eqref{Coupled_BDF2_1}-\eqref{Coupled_BDF2_3} and \eqref{gt3}-\eqref{gt2}, respectively. Define the error terms as 
\begin{align}
e_{\bm{u},f}^{n,i} := \bm{u}_{h,f}^{n,i} - \bm{u}_{h}^{n}, \quad e_{\xi,f}^{n,i} := \xi_{h,f}^{n,i} - \xi_{h}^{n}, \quad e_{p,f}^{n,i} := p_{h,f}^{n,i} - p_{h}^{n}. \label{errorf}
\end{align}
\cref{algo:GTIDA_bdf_fom} is globally convergent. 
There holds
    \begin{align}
            \sum_{n=2}^N \| \mathbb{D}_{t}^2 e_{\xi,f}^{n,i} \|_{L^2}^2  \leq K^2   \sum_{n=2}^N \| \mathbb{D}_{t}^2 e_{\xi,f}^{n,i-1} \|_{L^2}^2 ,
    \end{align}
     where the constant $K$ is given by
    \begin{align}
        K = \frac{1}{( \frac{c_0 \lambda }{\alpha^2} + 1) \times ( C(\tilde{\beta}, \mu)\lambda + 1 )}. \label{Kdef}
    \end{align}
    Since $K<1$, the iterative process is contractive.
\end{theorem}

\begin{proof}
We subtract \eqref{Coupled_BDF2_3} from \eqref{gt3}, then take $q_h = \mathbb{D}_{t}^2 e_{p,f}^{n,i}$ to obtain
\begin{align}
& a_3 ( \mathbb{D}_{t}^2 e_{p,f}^{n,i} , \mathbb{D}_{t}^2 e_{p,f}^{n,i} ) 
+ d ( e_{p,f}^{n,i}, \mathbb{D}_{t}^2 e_{p,f}^{n,i} )  = c( \mathbb{D}_{t}^2 e_{p,f}^{n,i} , \mathbb{D}_{t}^2 e_{\xi,f}^{n,i-1} ). \label{thm31:eq1}
\end{align}
With the notation of the bilinear forms, we can apply the Cauchy-Schwarz and Young's inequalities to get
\begin{align}
    & (c_0 + \frac{\alpha^2}{\lambda}) \| \mathbb{D}_{t}^2 e_{p,f}^{n,i} \|_{L^2}^2 + d ( e_{p,f}^{n,i}, \mathbb{D}_{t}^2 e_{p,f}^{n,i} )
    \nonumber \\
    & \leq \frac{1}{2} (c_0 + \frac{\alpha^2}{\lambda}) \| \mathbb{D}_{t}^2 e_{p,f}^{n,i} \|_{L^2}^2 + \frac{1}{2(c_0 + \frac{\alpha^2}{\lambda})} \frac{\alpha^2}{\lambda^2} \| \mathbb{D}_{t}^2 e_{\xi,f}^{n,i-1} \|_{L^2}^2. \label{thm31:eq2}
\end{align}
Recall the equality \eqref{bdfineq} to see that
\begin{align}
d ( e_{p,f}^{n,i}, \mathbb{D}_{t}^2 e_{p,f}^{n,i} ) = & \frac{k_p}{2\Delta t} \| [\nabla e_{p,f}^{n,i}; \nabla  e_{p,f}^{n-1,i}] \|_G^2 
- \frac{k_p}{2\Delta t} \| [\nabla  e_{p,f}^{n-1,i}; \nabla  e_{p,f}^{n-2,i}] \|_G^2 
\nonumber \\
& + \frac{k_p}{4\Delta t} \| \nabla  e_{p,f}^{n,i} - 2 \nabla  e_{p,f}^{n-1,i} + \nabla  e_{p,f}^{n-2,i} \|_{L^2}^2. \label{thm31:eq3}
\end{align}
Combine \eqref{thm31:eq2} and  \eqref{thm31:eq3}, we can deduce that
\begin{align}
    & (c_0 + \frac{\alpha^2}{\lambda}) \| \mathbb{D}_{t}^2 e_{p,f}^{n,i} \|_{L^2}^2 + \frac{k_p}{\Delta t} \| [\nabla e_{p,f}^{n,i}; \nabla  e_{p,f}^{n-1,i}] \|_G^2 
    \nonumber \\
    & \leq \frac{1}{(c_0 + \frac{\alpha^2}{\lambda})} \frac{\alpha^2}{\lambda^2} \| \mathbb{D}_{t}^2 e_{\xi,f}^{n,i-1} \|_{L^2}^2 + \frac{k_p}{\Delta t} \| [\nabla  e_{p,f}^{n-1,i}; \nabla  e_{p,f}^{n-2,i}] \|_G^2 . \label{thm31:eq4}
\end{align}
Applying the summation operator $\sum_{n=2}^N$ to \eqref{thm31:eq4}, with $\| [\nabla  e_{p,f}^{1,i}; \nabla  e_{p,f}^{0,i}] \|_G^2 = 0$, gives
\begin{align}
    & \sum_{n=2}^N  ( c_0 + \frac{\alpha^2}{\lambda} ) \| \mathbb{D}_{t}^2 e_{p,f}^{n,i} \|_{L^2}^2 + \frac{k_p}{\Delta t} \| [\nabla e_{p,f}^{N,i}; \nabla  e_{p,f}^{N-1,i}] \|_G^2
    \nonumber \\
    & \leq \sum_{n=2}^N \frac{1}{(c_0 + \frac{\alpha^2}{\lambda})} \frac{\alpha^2}{\lambda^2} \| \mathbb{D}_{t}^2 e_{\xi,f}^{n,i-1} \|_{L^2}^2. \label{thm31:eq5}
\end{align}
Discard the second positive term on the left-hand side of \eqref{thm31:eq5} to obtain
\begin{align}
    \sum_{n=2}^N \| \mathbb{D}_{t}^2 e_{p,f}^{n,i} \|_{L^2}^2 
    \leq \sum_{n=2}^N \Big( \frac{\alpha}{ c_0 \lambda 
    + \alpha^2} \Big)^2 \| \mathbb{D}_{t}^2 e_{\xi,f}^{n,i-1} \|_{L^2}^2. 
    \label{thm31:eq5s}
\end{align}

On the other hand, we subtract \eqref{Coupled_BDF2_1}, \eqref{Coupled_BDF2_2} from \eqref{gt1}, \eqref{gt2}, respectively, then apply the operator $\mathbb{D}_{t}^2$ to get
\begin{align}
& a_1( \mathbb{D}_{t}^2 e_{\bm{u},f}^{n,i},\bm{v}_h)-b(\bm{v}_h, \mathbb{D}_{t}^2 e_{\xi,f}^{n,i}) = 0, \label{thm31:eq6}
\\
& b( \mathbb{D}_{t}^2 e_{\bm{u},f}^{n,i},w_h) + a_2( \mathbb{D}_{t}^2 e_{\xi,f}^{n,i},w_h) = c(\mathbb{D}_{t}^2 e_{p,f}^{n,i},w_h),  \label{thm31:eq7}
\end{align}
After taking $\bm{v}_h = \mathbb{D}_{t}^2 e_{\bm{u},f}^{n,i}$ in \eqref{thm31:eq6} and $w_h = \mathbb{D}_{t}^2 e_{\xi,f}^{n,i}$ \eqref{thm31:eq7}, we sum them up to obtain
\begin{align}
a_1( \mathbb{D}_{t}^2 e_{\bm{u},f}^{n,i},\mathbb{D}_{t}^2 e_{\bm{u},f}^{n,i}) + a_2( \mathbb{D}_{t}^2 e_{\xi,f}^{n,i}, \mathbb{D}_{t}^2 e_{\xi,f}^{n,i})
= c(\mathbb{D}_{t}^2 e_{p,f}^{n,i}, \mathbb{D}_{t}^2 e_{\xi,f}^{n,i}),  \label{thm31:eq8}
\end{align}
Applying the discrete inf-sup condition \eqref{dinfsupcon} to \eqref{thm31:eq6} yields
\begin{align}
    \tilde{\beta} \| \mathbb{D}_{t}^2 e_{\xi,f}^{n,i} \|_{L^2} & \leq \sup_{\bm{v}_h \in \bm{V}_h} \frac{  b(\bm{v}_h, \mathbb{D}_{t}^2 e_{\xi,f}^{n,i} ) }{\| \bm{v}_h \|_{H^1}} \nonumber \\
    & = \sup_{\bm{v}_h \in \bm{V}_h} \frac{  a_1( \mathbb{D}_{t}^2 e_{\bm{u},f}^{n,i},\bm{v}_h) }{\| \bm{v}_h \|_{H^1}} \leq 2\mu C \| \varepsilon(\mathbb{D}_{t}^2 e_{\bm{u},f}^{n,i}) \|_{L^2}, \label{thm31:eq9}
\end{align}
which leads to the following inequality.
\begin{align}
    C(\tilde{\beta}, \mu) \| \mathbb{D}_{t}^2 e_{\xi,f}^{n,i} \|_{L^2}^2 \leq 2 \mu \| \varepsilon(\mathbb{D}_{t}^2 e_{\bm{u},f}^{n,i}) \|_{L^2}^2. \label{thm31:eq10}
\end{align}
We then apply the Cauchy inequality to \eqref{thm31:eq8}, together with \eqref{thm31:eq10}, to see that
\begin{align}
     ( C(\tilde{\beta}, \mu) + \frac{1}{\lambda} ) \| \mathbb{D}_{t}^2 e_{\xi,f}^{n,i} \|_{L^2}^2 \leq \frac{\alpha}{\lambda} \| \mathbb{D}_{t}^2 e_{\xi,f}^{n,i} \|_{L^2} \| \mathbb{D}_{t}^2 e_{p,f}^{n,i} \|_{L^2}. \label{thm31:eq11}
\end{align}
Assuming $\| \mathbb{D}_{t}^2 e_{\xi,f}^{n,i} \|_{L^2} \not= 0$, , we divide both sides of \eqref{thm31:eq11} by $\| \mathbb{D}_{t}^2 e_{\xi,f}^{n,i} \|_{L^2}$, then square both sides and sum over $n=2$ to $N$, which yields
\begin{align}
    \sum_{n=2}^N ( \lambda C(\tilde{\beta}, \mu) + 1 )^2 \| \mathbb{D}_{t}^2 e_{\xi,f}^{n,i} \|_{L^2}^2 \leq \sum_{n=2}^N \alpha^2 \| \mathbb{D}_{t}^2 e_{p,f}^{n,i} \|_{L^2}^2. \label{thm31:eq12}
\end{align}

At last, combine the two inequalities \eqref{thm31:eq5s} and \eqref{thm31:eq12}, we can find that
\begin{align}
    ( \lambda C(\tilde{\beta}, \mu) + 1 )^2 \sum_{n=2}^N \| \mathbb{D}_{t}^2 e_{\xi,f}^{n,i} \|_{L^2(\Omega)}^2 \leq \Big( \frac{\alpha^2}{ c_0 \lambda 
    + \alpha^2} \Big)^2 \sum_{n=2}^N \| \mathbb{D}_{t}^2 e_{\xi,f}^{n,i-1} \|_{L^2(\Omega)}^2.
\end{align}
This completes the proof.
\end{proof}

\begin{remark}[Discussion]
\label{remark1}
Together with \eqref{thm31:eq5s} and \eqref{thm31:eq8}, the above analysis implies that the relative error terms $e_{\bm{u},f}^{n,i}$, $e_{\xi,f}^{n,i}$, and $e_{p,f}^{n,i}$ defined in \eqref{errorf} converge to zero as $i \rightarrow \infty$.  Applying the triangular inequality, we see that the discrepancy between the exact solution and the numerical solution $(\bm{u}_h^n, \xi_h^n, p_h^n)$ generated by \cref{algo:Coupled_BDF2} serves as an upper bound for the absolute error between the exact solution and the iterative solution produced by \cref{algo:GTIDA_bdf_fom}.
\end{remark}

\section{POD-based reduced order modeling}
\label{Sec4}
In this section, we present a POD-based reduced order model aimed at reducing the computational cost of \cref{algo:GTIDA_bdf_fom}. The approach involves solving the generalized Stokes problem over a subset of time steps to generate training data for the ROM. Using these training data, a POD basis is constructed, prioritizing modes associated with the largest eigenvalues. The resulting POD basis is then employed to efficiently solve the entire generalized Stokes problem. Additionally, an error estimate analysis for the proposed POD-based algorithm is provided.

\subsection{Proper orthogonal decomposition}
\label{pod}

Let $\bm{\phi}_h \in \bm{V}_h$ represent the displacement solution field and $\rho_h \in W_h$ represent the total pressure solution field. Introduce an index set $\bm{\Lambda} = \{ n \ | \ \Lambda_n \in \{2, \cdots, N\} \} $, and $\#(\bm{\Lambda}) $ denotes the size of the index set $\bm{\Lambda}$. The sequences of displacement and pressure solutions are denoted by $\{\phi_h^{\Lambda_n}\}_{n \in \bm{\Lambda}}$ and $\{\rho_h^{\Lambda_n}\}_{n \in \bm{\Lambda}}$, respectively. 
These sequences serve as the snapshots for the Proper Orthogonal Decomposition (POD) method \cite{berkooz1993proper, hesthaven2016certified, ballarin2024projection}, which is used to construct the reduced basis (RB) by solving the following eigenvalue problems first: For $\textbf{C}^{\bm{\phi}_h}, \textbf{C}^{\rho_h} \in \mathbb{R}^{\#(\bm{\Lambda}) \times \#(\bm{\Lambda})}$ determine  $\gamma^{\bm{\phi}_h}, \gamma^{\rho_h} \in \mathbb{C}$ and $\underline{\textbf{v}}, \underline{\text{w}} \in \mathbb{R}^{\#(\bm{\Lambda})} \backslash \bm{0}$ such that
\begin{align*}
    \textbf{C}^{\bm{\phi}_h} \underline{\textbf{v}} = \gamma^{\bm{\phi}_h} \underline{\textbf{v}}, \quad \textbf{C}^{\rho_h} \underline{\text{w}} = \gamma^{\rho_h} \underline{\text{w}}, 
\end{align*}
where $[\textbf{C}^{\bm{\phi}_h}]_{nm} = (\bm{\phi}_h^{\Lambda_n},\bm{\phi}_h^{\Lambda_m})_1$, $[\textbf{C}^{\rho_h}]_{nm} = (\rho_h^{\Lambda_n},\rho_h^{\Lambda_m})$ for every $n, m \in \bm{\Lambda}$.
By construction, the matrices $\textbf{C}^{\bm{\phi}_h}$ and $\textbf{C}^{\rho_h}$ are symmetric positive definite, ensuring that their eigenvalues $\gamma^{\bm{\phi}_h}_1, \cdots, \gamma^{\bm{\phi}_h}_{\#(\bm{\Lambda})}$ and $\gamma^{\rho_h}_1, \cdots, \gamma^{\rho_h}_{\#(\bm{\Lambda})}$ are real and positive. Without loss of generality, we assume that the eigenvalues are sorted in a decreasing order. 
The eigenvectors associated with the eigenvalues $\gamma^{\bm{\phi}_h}_n$ and $\gamma^{\rho_h}_n$ are denoted as $\underline{\textbf{v}}_n$ and $\underline{\text{w}}_n$, respectively. 
The $n$-th POD modes are then defined as:
\begin{align*}
    \bm{\varphi}_n = \frac{1}{\sqrt{\gamma^{\bm{\phi}_h}_n}} \sum_{\beta = 1}^{{\#(\bm{\Lambda})}} [ {\underline{\textbf{v}}}_{n} ]_{\beta} \phi_h^{\Lambda_\beta}, \quad
    \varrho_n = \frac{1}{\sqrt{\gamma^{\rho_h}_n}} \sum_{\beta = 1}^{{\#(\bm{\Lambda})}} [ {\underline{\text{w}}}_{n} ]_{\beta} \rho_h^{\Lambda_\beta}.
\end{align*}
The resulting POD modes $\{\bm{\varphi}_n\}_{n \in \bm{\Lambda}}$ are orthonormal with respect to the $(\cdot, \cdot)_1$ inner product, while $\{\varrho_n\}_{n \in \bm{\Lambda}}$  are orthonormal with respect to the $(\cdot, \cdot)$ inner product.

To construct the RB spaces, we truncate the modes by retaining only the first $N_r$ displacement modes $\bm{\varphi}_1, \cdots, \bm{\varphi}_{N_r}$ and the first $N_r$ pressure modes $\varrho_1, \cdots, \varrho_{N_r}$. This truncation is typically guided by the decay of the eigenvalues, where the retained modes correspond to the largest eigenvalues, capturing the most significant features of the snapshots. The resulting truncated RB spaces are then defined as:
\begin{align*}
    \bm{V}_{h,r} = \text{span} ( \bm{\varphi}_1, \cdots, \bm{\varphi}_{N_r} ), \quad W_{h,r} =  \text{span}( \varrho_1, \cdots, \varrho_{N_r} ).
\end{align*}
These reduced spaces are then employed to project the governing equations of the full order model onto lower-dimensional spaces, achieving a significant reduction in computational cost while preserving the essential dynamics and accuracy of the numerical solutions.

\subsection{Fully discrete finite element methods of reduced order model}
\label{sec42}

In this section, we introduce a reduced order modeling (ROM) strategy to reduce the computational cost of \cref{algo:GTIDA_bdf_fom} via the Proper Orthogonal Decomposition (POD) method outlined in \cref{pod}. The ROM specifically targets the Generalized Stokes problem \eqref{gt1}–\eqref{gt2}, which is computationally intensive in its full order form. At the $i$-th iteration, snapshots $\{(\tilde{\bm{u}}_{h}^{\Lambda_n,i},\tilde{\xi}_{h}^{\Lambda_n,i})\}_{n \in \bm{\Lambda}} \in \bm{V}_h \times W_h$ are generated by solving the full order model (FOM). 
The reduced basis (RB) spaces are constructed via POD:
\begin{align}
    \bm{V}_{h,r}^i = \text{POD}(\{\tilde{\bm{u}}_{h}^{\Lambda_n,i}\}_{n \in \bm{\Lambda}}; N_r), \quad W_{h,r}^i = \text{POD}(\{\tilde{\xi}_{h}^{\Lambda_n,i}\}_{n \in \bm{\Lambda}};N_r),
\end{align}
where $\text{POD}(\cdot ; N_r)$ denotes the truncation to the first $N_r$ modes.
Assume that the corresponding eigenvalues are sorted in descending order as $\gamma^{\bm{\phi}_h,i}_1, \cdots, \gamma^{\bm{\phi}_h,i}_{\#(\bm{\Lambda})}$ and $\gamma^{\rho_h,i}_1, \cdots, \gamma^{\rho_h,i}_{\#(\bm{\Lambda})}$.
Rather than seeking solutions $(\bm{u}_{h,f}^{n,i},\xi_{h,f}^{n,i})\in \bm{V}_{h}^i \times W_{h}$ as in the full order scheme (\cref{algo:GTIDA_bdf_fom}), we now look for solutions $(\bm{u}_{h,r}^{n,i},\xi_{h,r}^{n,i}) \in \bm{V}_{h,r}^i \times W_{h,r}^i$ within the reduced spaces for $n = 2, \cdots, N$. This strategy ensures a substantial reduction in computational cost while maintaining accuracy. The proposed ROM-based global-in-time iterative decoupled BDF-$2$ algorithm is summarized in \cref{algo:GTIDA_bdf_rom}.

\begin{algorithm}
\raggedright
  \caption{: A ROM-based global-in-time iterative decoupled BDF-$2$ algorithm}
  \label{algo:GTIDA_bdf_rom}
  \textbf{Input:} initial data $\{\xi_{h,r}^{n}\}_{n=0}^{1} \subset W_h$ and $\{p_{h,r}^{n}\}_{n=0}^{1} \subset M_h$, initial guesses 
  \\
  \quad \quad \quad \ \ 
  $\{\xi_{h,r}^{n,0}\}_{n=2}^{N} \subset W_h$, initial iteration number $i = 0$, an index set $\bm{\Lambda}$,
  \\
  \quad \quad \quad \ \ and the reduced basis size of $N_r< \#(\bm{\Lambda}) $.
  \\
  \textbf{Output:} solutions $\{(\bm{u}_{h,r}^{n},\xi_{h,r}^{n},p_{h,r}^{n})\}_{n=2}^N \subset \bm{V}_h \times W_h \times M_h$.
  \\
  \textbf{while} not converged \textbf{do}
  \\
  \quad \quad set $i = i + 1$,  $p_{h,r}^{0,i} = p_h^0$, $p_{h,r}^{1,i} = p_{h,r}^1$, $\xi_{h,r}^{0,i-1} = \xi_{h,r}^{0}$, and $\xi_{h,r}^{1,i-1} = \xi_{h,r}^{1}$. 
  \\
  \quad \quad \textbf{Step a:} find $\{ p_{h,r}^{n,i} \}_{n=2}^{N} \subset M_h$ such that
  \vspace{-0.3cm}
  \begin{align}
    & a_3 ( \mathbb{D}_{t}^2 p_{h,r}^{n,i} , q_h ) 
    + d ( p_{h,r}^{n,i},q_h )  =  ( g^{n} , q_h ) + c( q_h , \mathbb{D}_{t}^2 \xi_{h,r}^{n,i-1} ), \quad \forall q_h \in M_h. \label{gtr1}
  \end{align}
  \vspace{-0.7cm}
  \\
  \quad \quad \textbf{Step b:} find snapshots $\{(\tilde{\bm{u}}_{h}^{\Lambda_n,i},\tilde{\xi}_{h}^{\Lambda_n,i})\}_{n \in \bm{\Lambda}} \in \bm{V}_h \times W_h$ such that
  \vspace{-0.3cm}
  \begin{align}
    & a_1(\tilde{\bm{u}}_{h}^{\Lambda_n,i},\bm{v}_h)-b(\bm{v}_h,\tilde{\xi}_{h}^{\Lambda_n,i}) = ( \bm{f}^{\Lambda_n},\bm{v}_h ), \quad \quad   \forall \bm{v}_h \in \bm{V}_h, \label{gtr2}
    \\
    & b(\tilde{\bm{u}}_{h}^{\Lambda_n,i},w_h) + a_2(\tilde{\xi}_{h}^{\Lambda_n,i},w_h) = c(p_{h,r}^{\Lambda_n,i},w_h)  ,\quad \forall w_h \in W_h. \label{gtr3}
  \end{align}
  \vspace{-0.7cm}
  \\
  \quad \quad \textbf{Step c: (ROM)} use the obtained snapshots $\{(\tilde{\bm{u}}_{h}^{\Lambda_n,i},\tilde{\xi}_{h}^{\Lambda_n,i})\}_{n \in \bm{\Lambda}}$ to construct the RB spaces  $\bm{V}_{h,r}^i$ and $W_{h,r}^i$ via the POD method, with a reduced dimension of size $N_r$. Then, compute all solutions $\{(\bm{u}_{h,r}^{n,i},\xi_{h,r}^{n,i})\}_{n=2}^N \subset \bm{V}_{h,r}^i \times W_{h,r}^i$ such that
  \vspace{-0.3cm}
  \begin{align}
    & a_1(\bm{u}_{h,r}^{n,i},\bm{v}_{h,r})-b(\bm{v}_{h,r},\xi_{h,r}^{n,i}) = ( \bm{f}^{n},\bm{v}_{h,r} ), \quad \quad \forall \bm{v}_{h,r} \in \bm{V}_{h,r}^i, \label{gtr4}
    \\
    & b(\bm{u}_{h,r}^{n,i},w_{h,r}) + a_2(\xi_{h,r}^{n,i},w_{h,r}) = c(p_{h,r}^{n,i},w_{h,r})  ,\quad \forall w_{h,r} \in W_{h,r}^i. \label{gtr5}
  \end{align}
  \vspace{-0.7cm}
  \\
  \textbf{end while}
  \\
  set $(\bm{u}_{h,r}^{n},\xi_{h,r}^{n},p_{h,r}^{n}) = (\bm{u}_{h,r}^{n,i},\xi_{h,r}^{n,i},p_{h,r}^{n,i})$ for $n = 2, \cdots, N$.
\end{algorithm}

\cref{algo:GTIDA_bdf_rom} can be regarded as equivalent to \cref{algo:GTIDA_bdf_fom} if the index set $\bm{\Lambda}$ includes all time steps, i.e., $\{\Lambda_n\}_{n\in \bm{\Lambda}} = \{2, 3, \cdots, N-1, N\}$, and the reduced basis dimension is set to $N_r=\#(\bm{\Lambda}) $, without any truncation. However, in practical applications, the index set $\bm{\Lambda}$ is usually selected to be much smaller, and truncation is applied to $N_r$, resulting in reduced order modeling (ROM) error. To evaluate the accuracy and effectiveness of the proposed ROM-based approach (\cref{algo:GTIDA_bdf_rom}), we analyze its error bounds as established in \Cref{Thm:rom}.

\begin{theorem}
\label{Thm:rom}
Let $(\bm{u}_h^n,\xi_h^n,p_h^n)$ and $(\tilde{\bm{u}}_{h}^{n,i},\tilde{\xi}_{h}^{n,i},\bm{u}_{h,r}^{n,i},\xi_{h,r}^{n,i},p_{h,r}^{n,i})$ be the solutions of \eqref{Coupled_BDF2_1}-\eqref{Coupled_BDF2_3} and \eqref{gtr1}-\eqref{gtr5}, respectively.
Define the error terms as 
\begin{align}
e_{\bm{u},r}^{n,i} := \bm{u}_{h,r}^{n,i} - \bm{u}_{h}^{n}, \quad e_{\xi,r}^{n,i} := \xi_{h,r}^{n,i} - \xi_{h}^{n}, \quad e_{p,r}^{n,i} := p_{h,r}^{n,i} - p_{h}^{n}. 
\label{errorr}
\end{align}
Then, the error in $\xi_{h,r}^{n,i}$ is bounded as follows:
\begin{align}
    \sum_{n=2}^N \| \mathbb{D}_{t}^2 e_{\xi,r}^{n,i} \|_{L^2}^2 \leq  C_1 \epsilon_{ROM}^i(\bm{\Lambda},N_r) + \frac{1+K^2}{2} \sum_{n=2}^N \| \mathbb{D}_{t}^2 e_{\xi,r}^{n,i-1} \|_{L^2}^2, \label{thm:rom1}
\end{align}
where $\epsilon_{ROM}^i(\bm{\Lambda},N_r) := \max\limits_{2 \leq n \leq N} \| \xi_{h,r}^{n,i} - \tilde{\xi}_{h}^{n,i} \|_{L^2}^2$ measures the reduced order modeling error at the $i$-th iteration, $K$ is defined in \eqref{Kdef}. Assuming the existence of an upper bound
\begin{align} 
\epsilon_{ROM}(\bm{\Lambda}, N_r) := \max\limits_{i \geq 0} \epsilon_{ROM}^i(\bm{\Lambda}, N_r), 
\end{align}
it follows that:
\begin{align}
    \lim_{i \rightarrow \infty}  \sum_{n=2}^N \| \mathbb{D}_{t}^2 e_{\xi,r}^{n,i} \|_{L^2}^2 \leq C_2 \epsilon_{ROM}(\Lambda, N_r). \label{thm:rom2}
\end{align}
Here, $C_1$ and $C_2$ are constants depending on $K$, $N$ and $\Delta t$. The term $\epsilon_{ROM}^i(\bm{\Lambda},N_r)$ represents the asymptotic reduced order modeling error at the $i$-th iteration, which is determined by the selected ROM parameters $\bm{\Lambda}$ and the reduced basis size $N_r$.
\end{theorem}

\begin{proof}
First, we focus on the error term $e_{\xi,r}^{n,i}$, which can be decomposed as: 
\begin{align}
    e_{\xi,r}^{n,i} = ( \xi_{h,r}^{n,i} - \tilde{\xi}_{h}^{n,i} ) + ( \tilde{\xi}_{h}^{n,i} - \xi_{h}^{n} ). \label{thm2:eq1}
\end{align}
Applying the operator $\mathbb{D}_{t}^2$ to \eqref{thm2:eq1} results in the following inequality.
\begin{align}
    \| \mathbb{D}_{t}^2 e_{\xi,r}^{n,i} \|_{L^2} & \leq \| \mathbb{D}_{t}^2 (\xi_{h,r}^{n,i} - \tilde{\xi}_{h}^{n,i}) \|_{L^2} + \| \mathbb{D}_{t}^2 ( \tilde{\xi}_{h}^{n,i} - \xi_{h}^{n} ) \|_{L^2}. \label{thm2:eq2}
\end{align}
Squaring both sides of \eqref{thm2:eq2} and applying Young's inequality yields
\begin{align}
    \| \mathbb{D}_{t}^2 e_{\xi,r}^{n,i} \|_{L^2}^2 \leq (1+\frac{1}{\epsilon}) \| \mathbb{D}_{t}^2 (\xi_{h,r}^{n,i} - \tilde{\xi}_{h}^{n,i}) \|_{L^2}^2 + (1+\epsilon) \| \mathbb{D}_{t}^2 ( \tilde{\xi}_{h}^{n,i} - \xi_{h}^{n} ) \|_{L^2}^2, \label{thm2:eq3}
\end{align}
where $\epsilon>0$ is a parameter to be determined later.
To estimate the reduced order error term, we proceed as follows.
\begin{align}
    & \sum_{n=2}^N \| \mathbb{D}_{t}^2 (\xi_{h,r}^{n,i} - \tilde{\xi}_{h}^{n,i}) \|_{L^2}^2 \leq C \frac{1}{\Delta t^2} \sum_{n=2}^N \| \xi_{h,r}^{n,i} - \tilde{\xi}_{h}^{n,i} \|_{L^2}^2 
    \nonumber \\
    & \leq C \frac{N}{\Delta t^2} \max\limits_{2 \leq n \leq N} \| \xi_{h,r}^{n,i} - \tilde{\xi}_{h}^{n,i} \|_{L^2}^2 = C \frac{N \epsilon_{ROM}^i(\bm{\Lambda},N_r) }{\Delta t^2}. \label{thm2:eq5}
\end{align}
By employing \cref{Thm1}, the following inequality holds.
\begin{align}
    \sum_{n=2}^N \| \mathbb{D}_{t}^2 ( \tilde{\xi}_{h}^{n,i} - \xi_{h}^{n} ) \|_{L^2}^2  \leq K^2   \sum_{n=2}^N \| \mathbb{D}_{t}^2 ( \xi_{h,r}^{n,i-1} - \xi_{h}^{n} ) \|_{L^2}^2.  \label{thm2:eq4}
\end{align}
Next, applying the summation operator $\sum_{n=2}^N$ to \eqref{thm2:eq3}, together with \eqref{thm2:eq5} and \eqref{thm2:eq4}, yields
\begin{align}
    \sum_{n=2}^N  \| \mathbb{D}_{t}^2 e_{\xi,r}^{n,i} \|_{L^2}^2 \leq C (1+\frac{1}{\epsilon}) \frac{N \epsilon_{ROM}^i(\bm{\Lambda},N_r) }{\Delta t^2} + (1+\epsilon) K^2 \sum_{n=2}^N  \| \mathbb{D}_{t}^2 e_{\xi,r}^{n,i-1} \|_{L^2}^2. \label{thm2:eq6}
\end{align}
Finally, by choosing $\epsilon = \frac{1}{2K^2} - \frac{1}{2}$, the following inequality is obtained:
\begin{align}
    \sum_{n=2}^N \| \mathbb{D}_{t}^2 e_{\xi,r}^{n,i} \|_{L^2}^2 \leq  C_1 \epsilon_{ROM}^i(\bm{\Lambda},N_r) + \frac{1+K^2}{2} \sum_{n=2}^N \| \mathbb{D}_{t}^2 e_{\xi,r}^{n,i-1} \|_{L^2}^2, \label{thm2:eq7}
\end{align}
where $C_1 = C(K,N,\Delta t)$ is a constant. This provides the first result \eqref{thm:rom1}. Since $K<1$, it follows that $\frac{1+K^2}{2}<1$. Assuming the existence of an upper bound $ \epsilon_{ROM}(\bm{\Lambda},N_r) \geq \epsilon_{ROM}^i(\bm{\Lambda},N_r)$ for all $i\geq1$, a recursive argument yields the following estimate:
\begin{align}
    \sum_{n=2}^N \| \mathbb{D}_{t}^2 e_{\xi,r}^{n,i} \|_{L^2}^2 \leq  C_2 \epsilon_{ROM}(\bm{\Lambda},N_r) + \left( \frac{1+K^2}{2} \right)^i \sum_{n=2}^N \| \mathbb{D}_{t}^2 e_{\xi,r}^{n,0} \|_{L^2}^2, \label{thm2:eq8}
\end{align}
where $C_2 = 2C_1/(1-K^2)$ is a constant. This leads to the second result \eqref{thm:rom2}.
\end{proof}

\begin{remark}
    The index set $\bm{\Lambda} = \{ n \ | \ \Lambda_n \in \{2, \cdots, N\} \}$ is typically chosen with $\Lambda_1 = 2$ and $\Lambda_{\#(\bm{\Lambda})} = N$ to ensure that the POD method captures the solution dynamics across the entire time interval. 
    In the BDF-$m$ scheme, the first $m$ initial solution vectors can be incorporated into the snapshot matrix to enrich the reduced basis, thereby enhancing the accuracy of the reduced order model. 
\end{remark}

\section{Numerical experiments}
\label{Sec5}
In this section, we present numerical experiments to evaluate the efficiency and accuracy of the proposed algorithms, aiming to verify the second-order temporal accuracy and analyze the influence of the choices for $\bm{\Lambda}$ and $N_r$ on the performance of the reduced order model. The implementation is carried out using the FEniCSx\footnote{FEniCSx: \url{https://github.com/FEniCS/dolfinx}}, RBniCSx\footnote{RBniCSx: \url{https://github.com/RBniCS/RBniCSx}} and Multiphenicsx\footnote{Multiphenicsx: \url{https://multiphenics.github.io}} libraries \cite{baratta2023dolfinx, RozzaBallarinScandurraPichi2024}. 
\subsection{Example 1}
\label{Ex1}

We consider a unit square domain $\Omega=[0,1]^2$ with the final time $T=1.0$. The body force $\bm{f}$, source term $g$, initial conditions, and mixed Neumann and Dirichlet boundary conditions are chosen so that the exact solution is:
\begin{align*}
    & \bm{u} = 
    \begin{bmatrix}
    \sin{(\pi x t)} \cos{(\pi y t)} \\
    \cos{(\pi x t)} \sin{(\pi y t)}
    \end{bmatrix} x y (1-x)^2(1-y), 
    \\
    & p = \cos{(t + x - y)} x y (1-x)^2(1-y),
\end{align*}
where we impose $\Gamma_\sigma = \Gamma_q =\{(1,y); 0\leq y\leq1\}$ and $\Gamma_u = \Gamma_p = \partial \Omega \setminus \Gamma_\sigma$.
The physical parameters are specified as follows:
\begin{align*}
    \lambda=10^2, \quad \mu = 10^2, \quad \alpha = 1.0, \quad c_0=10^{-2}, \quad k_p=10^{-2}.
\end{align*}
The computational mesh is uniform with $h = \frac{1}{32}$, and we use the finite element spaces $(\bm{P}_4, P_3, P_4)$ for displacement, total pressure, and pressure, respectively.
Consistent with \cref{Thm1} and \cref{Thm:rom}, our analysis focuses on two key error terms: (1) the absolute error of the numerical solutions obtained using \cref{algo:Coupled_BDF2}, and (2) the relative error between the solutions computed by \cref{algo:Coupled_BDF2} and those from \cref{algo:GTIDA_bdf_fom} (or \cref{algo:GTIDA_bdf_rom}).

To verify the temporal accuracy of \cref{algo:Coupled_BDF2}, we compute the $L^2$ error for the total pressure, together with the $H^1$ semi-errors for displacement and pressure, with various time step sizes ($\Delta t = \frac{1}{8} , \frac{1}{16} , \frac{1}{32} , \frac{1}{64}$). These results are presented in \cref{E1}. We observe that these errors decrease at a rate of $2$, which aligns with the optimal convergence rate expected when using the BDF-$2$ scheme. This result demonstrates the second-order temporal accuracy of \cref{algo:Coupled_BDF2}.

\begin{table}[ht]
\setlength{\belowcaptionskip}{0.1cm}
\begin{center}
\renewcommand\arraystretch{1.5}
\caption{Errors at final time $t=T$ and convergence rates of Algorithm \ref{algo:Coupled_BDF2} for Example 1.}
\label{E1}
\centering
{  \footnotesize
	\begin{tabular}{c|cc|cc|cc}
	\hline
	$\Delta t$ & \multicolumn{1}{l}{$|\bm{u}_{h}^{N} - \bm{u}^N|_{H^1}$} & Orders & \multicolumn{1}{l}{$\|\xi_{h}^{N} - \xi^N\|_{L^2}$} & Orders & \multicolumn{1}{l}{$| p_{h}^{N} - p^N |_{H^1}$} & \multicolumn{1}{l}{Orders} \\ \hline
    1/8  & 5.354e-05 &      & 1.228e-02 &      & 1.126e-01 &      \\
    1/16 & 1.346e-05 & 1.99 & 2.987e-03 & 2.04 & 2.820e-02 & 2.00 \\
    1/32 & 3.355e-06 & 2.00 & 7.304e-04 & 2.03 & 7.006e-03 & 2.01 \\
    1/64 & 8.428e-07 & 1.99 & 1.808e-04 & 2.01 & 1.744e-03 & 2.01 \\
    \hline
	\end{tabular}
}
\end{center}
\end{table}

\begin{figure}[htbp]
\centering
\subfigure{\includegraphics[width=0.78\textwidth]{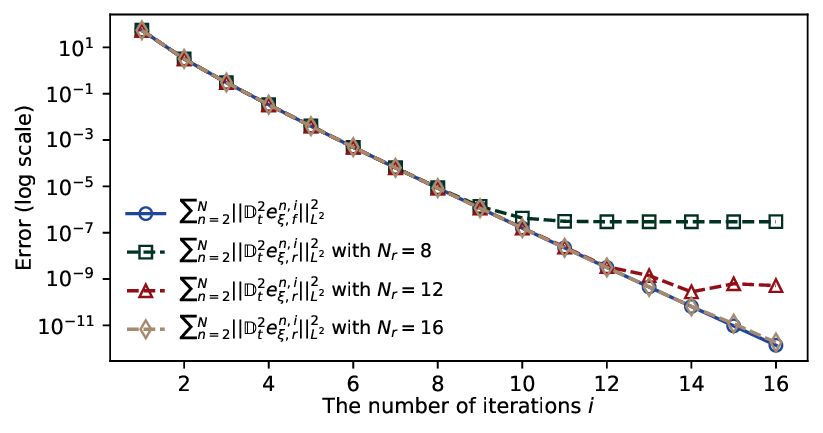}}
\caption{Example 1. Convergence behavior of the analyzed error terms for \cref{algo:GTIDA_bdf_fom} and \cref{algo:GTIDA_bdf_rom}, using the full index set $\bm{\Lambda}^f$ and reduced basis dimensions $N_r = 8, 12, 16$.}
\label{Fig1}
\end{figure}

After verifying the convergence rate of numerical result obtained from the coupled algorithm (\cref{algo:Coupled_BDF2}), we fix the time step size $\Delta t = 1/64$ and evaluate the performance of the iterative algorithms (\cref{algo:GTIDA_bdf_fom} and \cref{algo:GTIDA_bdf_rom}). An important objective is to validate the conclusions established in \cref{Thm1} and \cref{Thm:rom}. To achieve this, \cref{algo:GTIDA_bdf_rom} is constructed using the full index set $\bm{\Lambda}^f$, defined as $\Lambda_n^f = n+1$ for $n=1, \cdots, N-1$ (yielding 63 snapshots), while varying the reduced basis dimension ($N_r = 8, 12, 16$). Here, the full index set is specifically chosen to isolate and analyze the influence of $N_r$ without introducing additional variability caused by the selection of the index set. 
The results, summarized in \cref{Fig1}, confirm the convergence behavior of the error terms discussed above. Specifically, we monitor the quantity $\sum_{n=2}^N \| \mathbb{D}_t^2 e_{\xi,\star}^{n,i} \|_{L^2}^2$, where $e_{\xi,\star}^{n,i} = \xi_{h,\star}^{n,i} - \xi_h^n$ with $\star=f$ for \cref{algo:GTIDA_bdf_fom} and $\star=r$ for \cref{algo:GTIDA_bdf_rom}. 
Specifically, the figure shows the decay of the monotonically decreasing error term for the full order modeling approach (\cref{algo:GTIDA_bdf_fom}), which achieves the best performance, as expected. Meanwhile, the convergence behaviors of the corresponding error terms for the reduced order modeling approach (\cref{algo:GTIDA_bdf_rom}) with different reduced basis dimensions initially resemble those of \cref{algo:GTIDA_bdf_fom} during the early iterations. The errors for the reduced order modeling approach stabilize after a certain point, with larger values of $N_r$ resulting in smaller stabilized errors. This demonstrates that reduced order modeling achieves sufficiently accurate results.

\begin{figure}[htbp]
\centering
\subfigure{\includegraphics[width=0.49\textwidth]{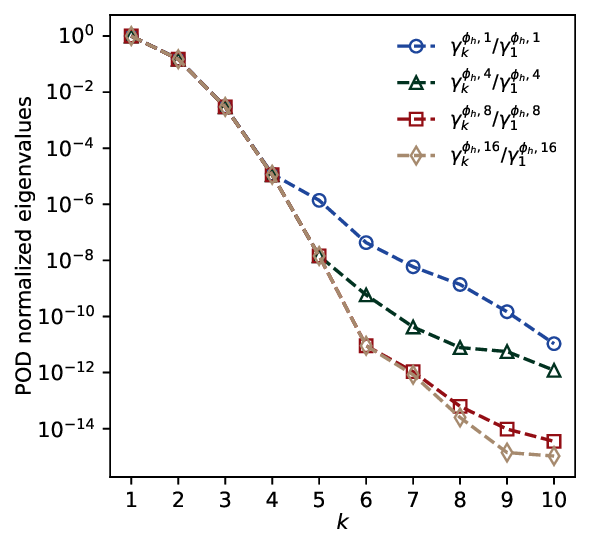}} 
\subfigure{\includegraphics[width=0.49\textwidth]{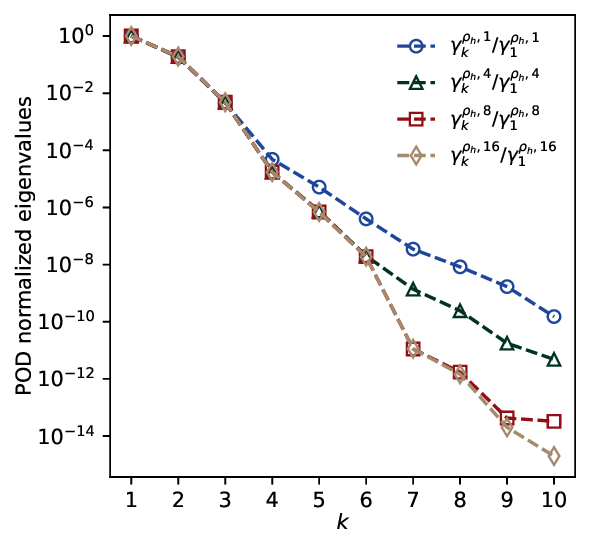}}
\caption{
Example 1. Normalized eigenvalues of the reduced basis spaces $\bm{V}_{h,r}^i$ (left) and $W_{h,r}^i$ (right) at iterations $i = 1, 4, 8, 16$, using the full index set $\bm{\Lambda}^f$ and reduced basis dimensions $N_r = 16$.}
\label{Fig1e}
\end{figure}

Moreover, for the case $N_r = 16$, we present the normalized eigenvalues of the reduced basis spaces $\bm{V}_{h,r}^i$ (left panel) and $W_{h,r}^i$ (right panel) at iterations $i = 1, 4, 8, 16$, as shown in Figure \ref{Fig1e}. 
As mentioned in Section \ref{sec42}, the ordered eigenvalues of $\bm{V}_{h,r}^i$ and $W_{h,r}^i$ are denoted as $\gamma^{\bm{\phi}_h,i}_k$ and $\gamma^{\rho_h,i}_k$ ($k=1,\cdots,\#(\Lambda)$), where $\gamma^{\bm{\phi}_h,i}_1$ and $\gamma^{\rho_h,i}_1$ represent the largest eigenvalues of the respective spaces.
The normalized eigenvalues, computed as $\gamma^{\bm{\phi}_h,i}_k / \gamma^{\bm{\phi}_h,i}_1$ and $\gamma^{\rho_h,i}_k / \gamma^{\rho_h,i}_1$, exhibit a clear decay pattern.
Additionally, the decay trend becomes convergent as the iteration number increases, indicating that the reduced basis spaces progressively stabilize.

\begin{figure}[htbp]
\centering
\subfigure{\includegraphics[width=0.78\textwidth]{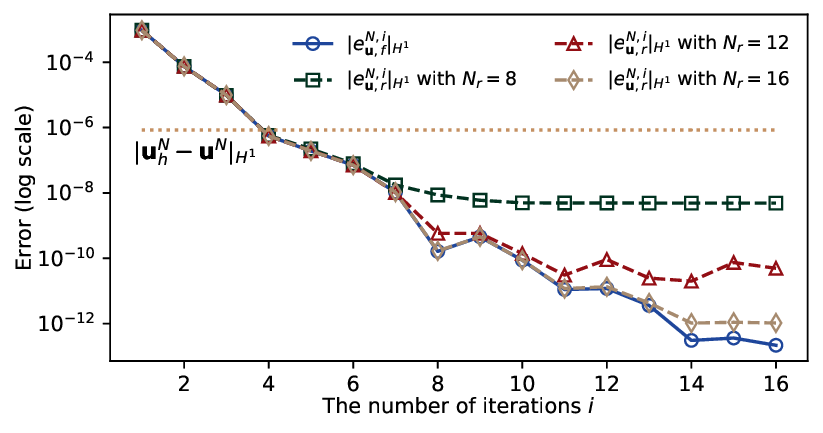}} 
\subfigure{\includegraphics[width=0.78\textwidth]{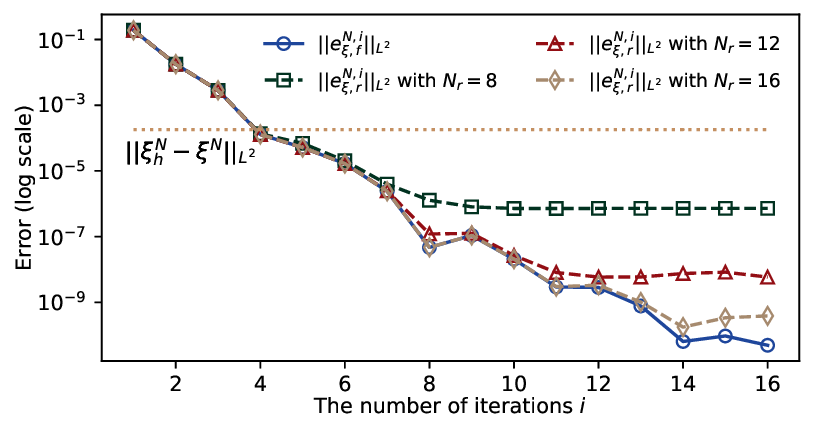}} 
\subfigure{\includegraphics[width=0.78\textwidth]{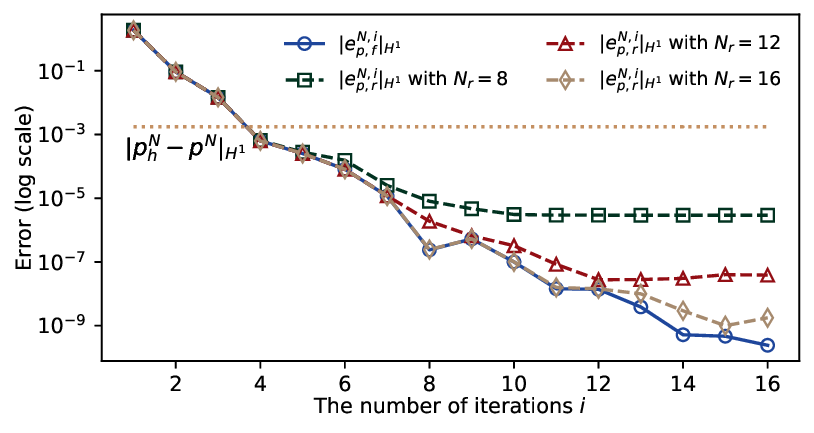}}
\caption{
Example 1. Convergence behaviors of the error terms at the final time for \cref{algo:GTIDA_bdf_fom} and \cref{algo:GTIDA_bdf_rom}, using the full index set $\bm{\Lambda}^f$ and reduced basis dimensions $N_r = 8, 12, 16$.}
\label{Fig2}
\end{figure}

Although \cref{Fig1} illustrates the convergence behavior of the analyzed error terms for the iterative algorithms, it remains challenging to assess whether the performance meets the desired level of accuracy. To provide a more detailed and clearer evaluation, we focus on the convergence behavior of the error terms specifically at the final time for \cref{algo:GTIDA_bdf_fom} and \cref{algo:GTIDA_bdf_rom}. 
The analysis presented in Figure~\ref{Fig2} enables a comprehensive comparison with the numerical error introduced by the finite element discretization, as reported in \cref{E1} for the case $\Delta t = \frac{1}{64}$ and depicted by the dotted lines in the following figure.

In \cref{Fig2}, the dashed lines show the evolution of the $L^2$ error term for total pressure and the $H^1$ semi-error terms for displacement and pressure under varying reduced basis dimensions ($N_r = 8, 12, 16$) with respect to the number of iterations. 
The solid lines confirm that the full order modeling approach (\cref{algo:GTIDA_bdf_fom}) consistently achieves the best performance.
For the reduced order modeling approach (\cref{algo:GTIDA_bdf_rom}), the convergence behavior closely follows the full order results during the initial iterations. 
However, as the iterations progress, the error terms stabilize, with larger reduced basis dimensions resulting in smaller error magnitudes. We also observe that tests with smaller reduced basis dimensions reach a ``stable'' state earlier, allowing the computations to be stopped sooner while still maintaining sufficient accuracy.
For a sufficient number of iterations, the iterative error for all methods falls below the discretization error, indicating that at some point, the finite element discretization error is dominating the overall error against the exact solution.
It highlights that the reduced order modeling strategy can achieve sufficiently accurate solutions with small reduced basis dimensions. 

\subsection{Example 2}
\label{Ex2}

\begin{figure}[ht]
\centering
\subfigure{
\includegraphics[width=0.450\textwidth,trim=0 0 0 0,clip]{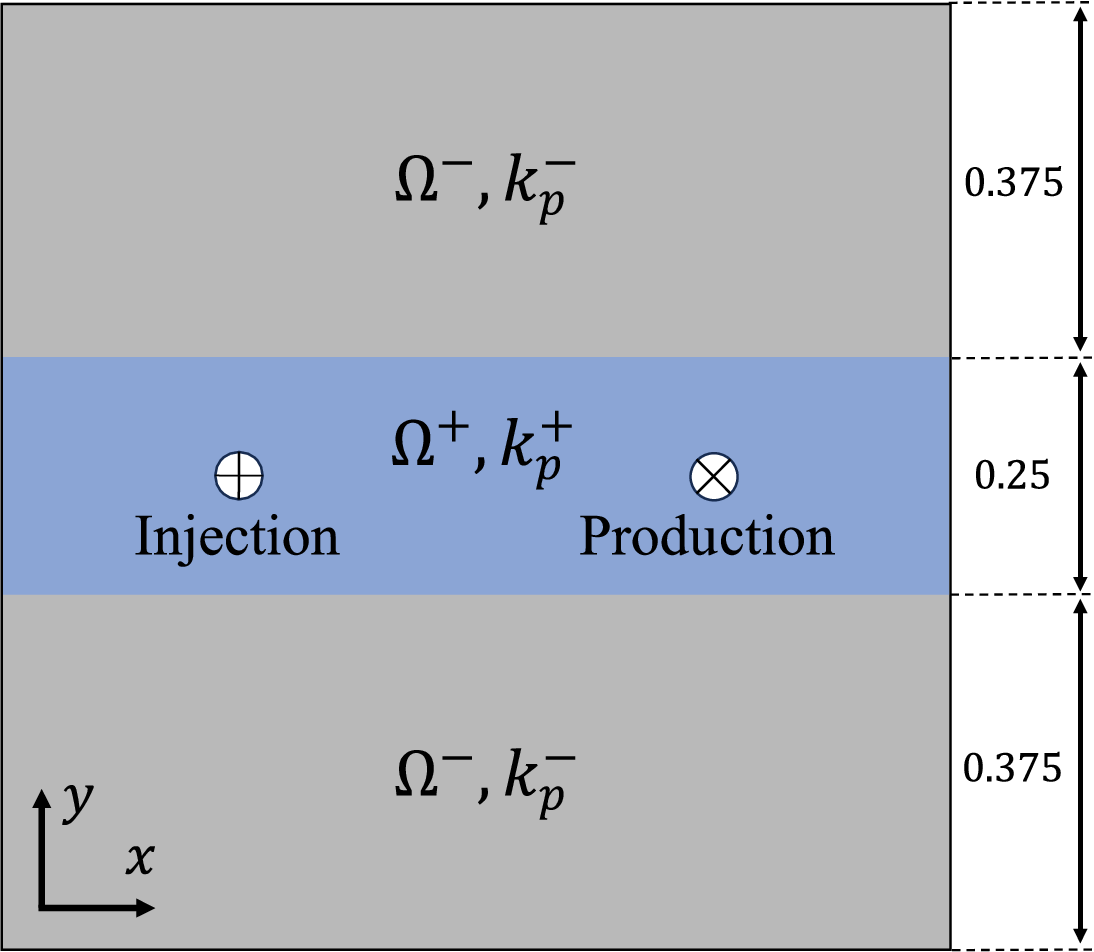}}
\caption{Example 2. 
Illustration of heterogeneous poroelastic medium with injection and production wells. 
}
\label{fig:EX2a}
\end{figure}

In this example, we study a problem involving injection and production processes in a heterogeneous poroelastic medium, following the setup described in \cite{yi2024physics, ballarin2024projection}. The computational domain is the unit square $\Omega = [0,1]^2$, featuring a vertically stratified permeability profile. Specifically, the domain is partitioned into two subregions:
\begin{align*} 
\Omega = \Omega^+ \cup \Omega^-, \quad \text{where} \quad 
\begin{cases} 
    \Omega^+ = \{(x,y) \in \Omega \mid 0.375 \leq y \leq 0.625\} \\ 
    \Omega^- = \Omega \setminus \Omega^+
\end{cases}.
\end{align*}
The middle subdomain $\Omega^+$ represents a high-permeability channel with $k_p^+ = 0.1$, while the surrounding subdomain $\Omega^-$ has a low permeability of $k_p^- = 10^{-5}$. We assume a vanishing body force, i.e., $\bm{f} = \bm{0}$. Injection and production wells are located at the points $(x_1, y_1) = (0.25, 0.5)$ and $(x_2, y_2) = (0.75, 0.5)$, respectively. The source term is defined as
\begin{equation*}
    g(x, y) = 10^{-2} \left( 
        e^{\left(-1000(x - x_1)^2 - 1000(y - y_1)^2\right)} 
        - e^{\left(-1000(x - x_2)^2 - 1000(y - y_2)^2\right)}
    \right),
\end{equation*}
modeling Gaussian injection and extraction profiles centered at the respective well locations.
An illustration of the domain and setup is shown in \cref{fig:EX2a}.

\begin{figure}[H]
\centering
\subfigure[$t = 0.0001$]{
\includegraphics[width=0.450\textwidth,trim=330 250 330 200,clip]{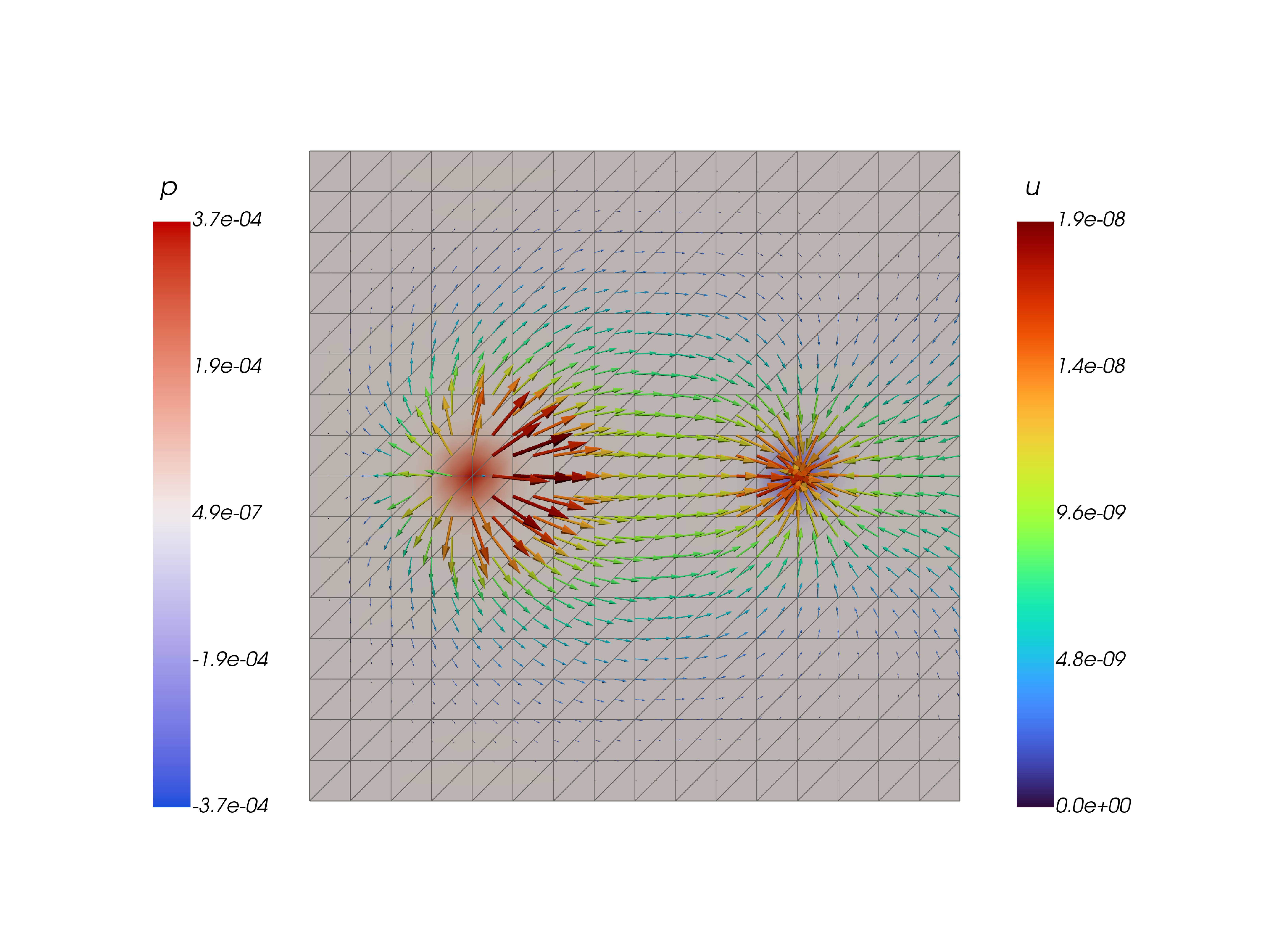}}
\subfigure[$t = 0.001$]{
\includegraphics[width=0.450\textwidth,trim=330 250 330 200,clip]{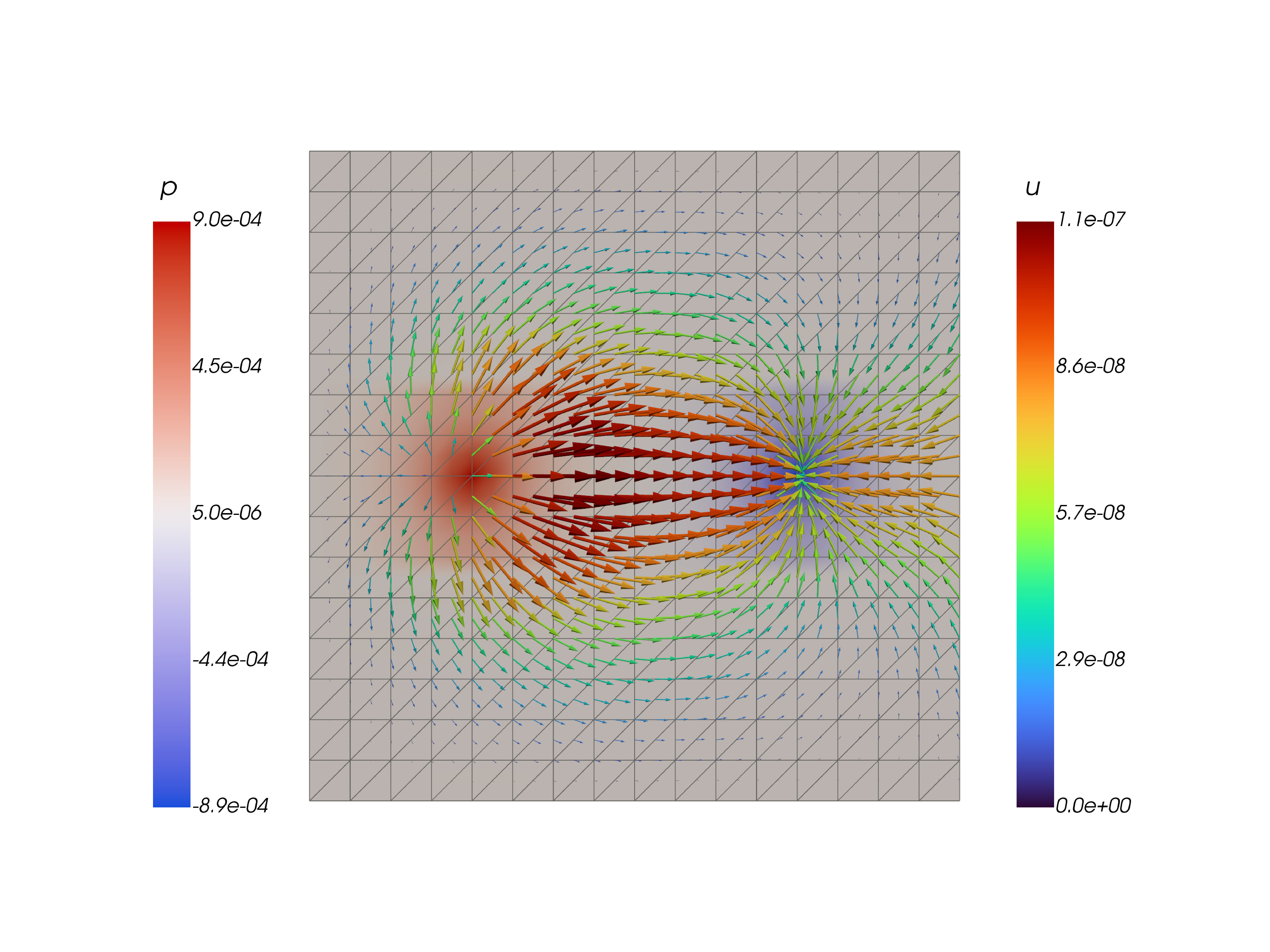}}
\subfigure[$t = 0.01$]{
\includegraphics[width=0.450\textwidth,trim=330 250 330 200,clip]{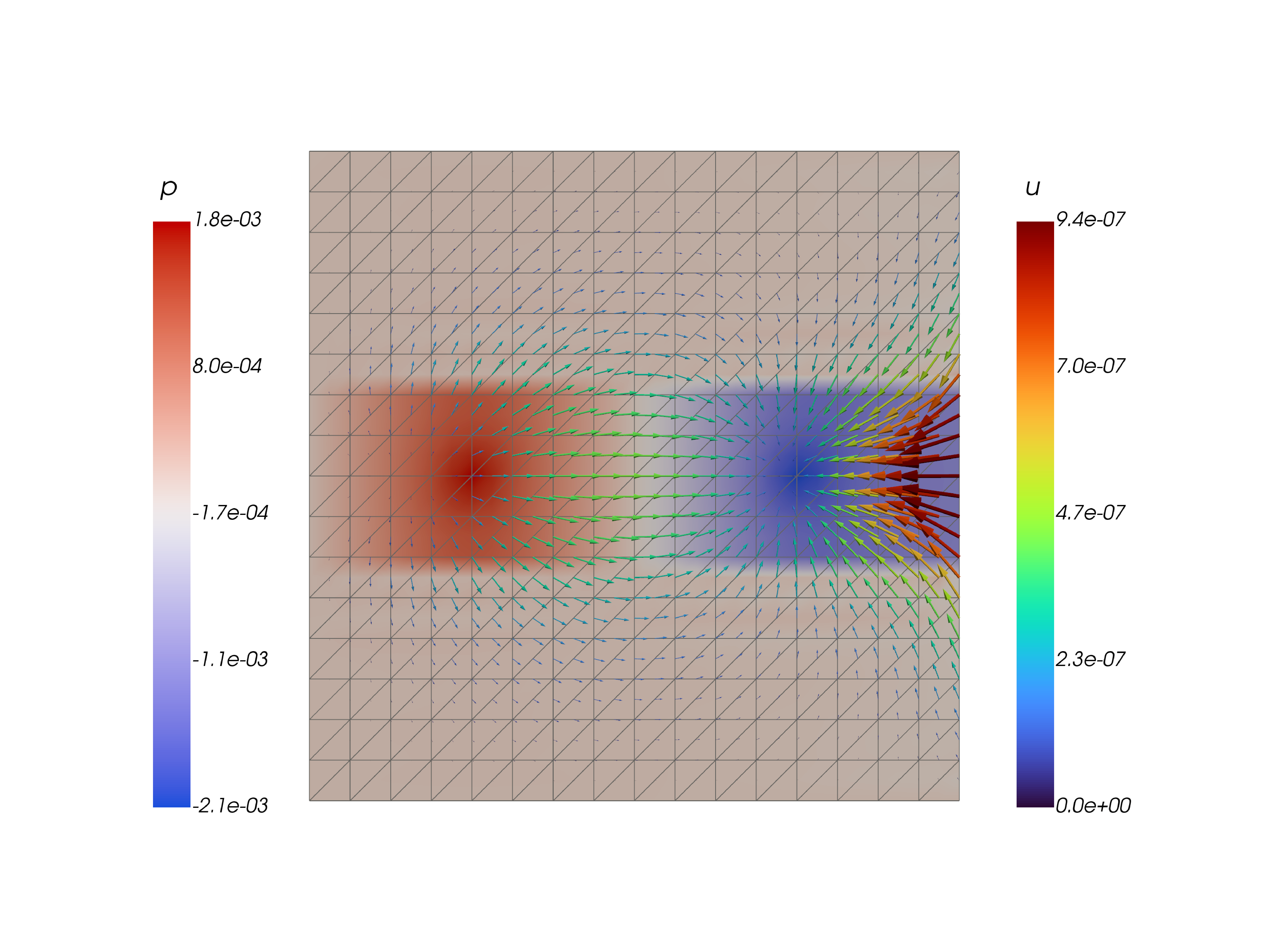}}
\caption{Example 2. Numerical results (displacement field $\bm{u}_h^N$ and pressure distribution $p_h^N$) obtained using \cref{algo:Coupled_BDF2} at different time snapshots: $(a)$$ t=0.0001$, $(b)$$ t=0.001$, $(c)$$ t=0.01$.
}
\label{fig:EX2}
\end{figure}

The physical parameters and boundary decomposition are identical to those used in \cref{Ex1}, with boundary conditions specified in \eqref{bc1}–\eqref{bc4}. Besides, the initial conditions for displacement and pressure are set to zero.
In this case, the spatial mesh is fixed with $h = \frac{1}{32}$, and the simulation is carried out up to final time $T = 0.01$ using time step $\Delta t = \frac{T}{256}$. The Taylor–Hood element $(\bm{P}_3, P_2)$ is used for the pair of displacement and total pressure, while the Lagrange element $P_3$ is used for the pressure. Since the BDF-$2$ method requires solution data at time step $t^1$, we employ \cref{algo:Coupled_BDF2} with the BDF-$1$ method to compute the initial step. This provides the necessary input for Algorithms~\ref{algo:Coupled_BDF2}, \ref{algo:GTIDA_bdf_fom}, and \ref{algo:GTIDA_bdf_rom}, while preserving second-order local accuracy. To illustrate the complexity of this benchmark, \cref{fig:EX2} presents the numerical results obtained using \cref{algo:Coupled_BDF2}, showing the significant evolution of the displacement field $\bm{u}_h$ and pressure distribution $p_h$ at different time snapshots.

\begin{figure}[H] 
\centering 
\subfigure{\includegraphics[width=0.9\textwidth,trim=0 5 0 10,clip]{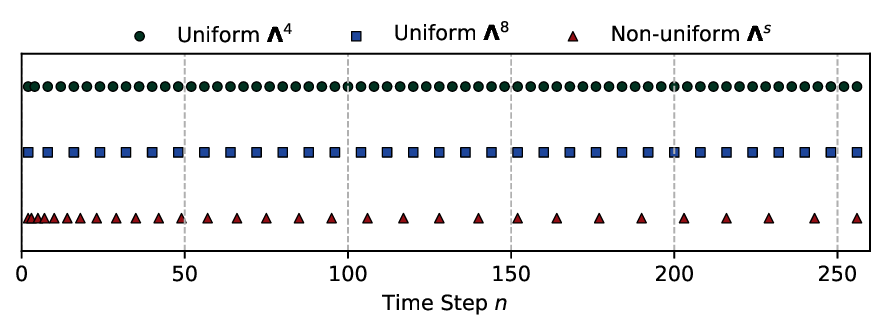}} 
\caption{Example 2. Comparison of sparse  index set selection strategies.} 
\label{ex2datapoint} 
\end{figure}

Since the primary objective of this study is to reduce computational cost, we have already investigated the impact of varying the reduced basis dimension. The remaining task is to examine whether using smaller index sets can maintain accuracy while further lowering the computational burden.  
We introduce two reduced index sets, $\bm{\Lambda}^4$ and $\bm{\Lambda}^8$, as alternatives to the full index set $\bm{\Lambda}^f$. These are defined as follows:
\begin{align*}
    \Lambda^4_n = \left\{
    \begin{array}{ll}
    2 & \quad n = 1 \\
    4({n-1}) & \quad 2 \leq n \leq \left(  \frac{N}{4}+1 \right)
    \end{array}
    \right.
    ,
    \quad
    \Lambda^8_n = \left\{
    \begin{array}{ll}
    2 & \quad n = 1 \\
    8({n-1}) & \quad 2 \leq n \leq \left( \frac{N}{8}+1 \right)
    \end{array}
    \right.
    .
\end{align*}
Since time-dependent poroelastic systems often exhibit rapid changes in state variables at early times and much slower variations at later times, it is beneficial to acquire more snapshots during the early stages and fewer snapshots later, as mentioned in \cite{siade2010snapshot}. 
Therefore, we also consider a non-uniform set $\bm{\Lambda}^s$ is constructed by mapping the first half of the Legendre polynomial roots to the integer interval from $2$ to $256$, followed by the removal of duplicate values.
Such approach allows for denser sampling in the early regions while maintaining a low snapshot budget. 
In total, four index sets are compared: the full index set $\bm{\Lambda}^f$ (255 snapshots), the uniform index sets $\bm{\Lambda}^4$ (65 snapshots) and $\bm{\Lambda}^8$ (33 snapshots), and the non-uniform set $\bm{\Lambda}^s$ (30 snapshots). The three sparse index sets are shown in \cref{ex2datapoint}. 

\begin{figure}[H]
\centering
\includegraphics[width=0.75\textwidth,trim=0 7 0 0,clip]{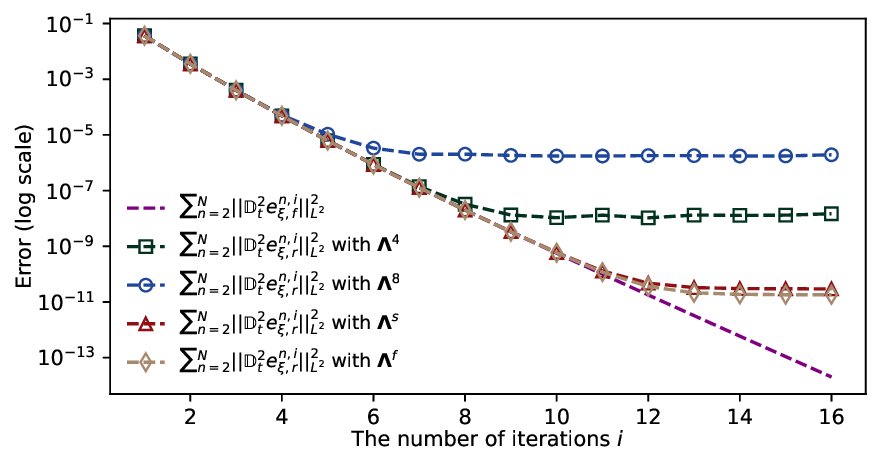}
\caption{Example 2. Convergence behavior for the analyzed error terms for \cref{algo:GTIDA_bdf_fom} and \cref{algo:GTIDA_bdf_rom}, comparing index sets $\bm{\Lambda}^4$, $\bm{\Lambda}^8$, $\bm{\Lambda}^s$, $\bm{\Lambda}^f$ with fixed reduced basis dimension $N_r = 16$. 
}
\label{figex2r}
\end{figure}

Finally, we carry out the numerical experiment to evaluate the performance of the proposed iterative algorithms. Here, we mainly focus on the summations $\sum_{n=2}^N \|\mathbb{D}_t^2 e_{\xi,\star}^{n,i}\|_{L^2}^2$ ($\star = f \ \text{or} \ r $) construed by varying index sets with fixed reduced basis dimension $N_r = 16$. The convergence behavior of these errors with respect to the iteration number $i$ is shown in \cref{figex2r}.

In \cref{figex2r}, the purple dashed line represents the error from the full order modeling approach (\cref{algo:GTIDA_bdf_fom}), which converges to zero as expected. Among the reduced order modeling results (\cref{algo:GTIDA_bdf_rom}), the full index set $\bm{\Lambda}^f$ serves as the reference benchmark, delivering the best performance. It is worth noting that the observed error plateau is due to the truncation at $N_r = 16$. Among the sparser index sets, the uniform set $\bm{\Lambda}^8$ yields the least accurate results, while $\bm{\Lambda}^4$, which uses roughly twice as many snapshots, demonstrates improved accuracy. Remarkably, the non-uniform set $\bm{\Lambda}^s$ performs even better than $\bm{\Lambda}^4$ and closely matches the accuracy of the full index set $\bm{\Lambda}^f$, despite using fewer snapshots than $\bm{\Lambda}^8$. These findings highlight that careful selection of the index set can significantly reduce the number of required snapshots while preserving high accuracy. 
It is also worth noting the computational cost. This test was executed in a cloud-based environment using Google Colab (Intel Xeon CPU @ 2.20GHz, 13GB RAM).
Each iteration of the full order algorithm (\cref{algo:GTIDA_bdf_fom}) took approximately 3 minutes to compute. By contrast, the reduced order algorithm with the full index set $\bm{\Lambda}^f$ took about 6 minutes per iteration, due primarily to the overhead introduced by the POD basis generation (step (c) in \cref{algo:GTIDA_bdf_rom}), which is computationally expensive. However, using a reduced index set $\bm{\Lambda}^s$, iteration time dropped to approximately 1 minute, achieving a threefold speedup compared to the full order approach. This gain is attributed to fewer snapshots and faster POD basis construction.
Furthermore, the similarity in error trends during the initial iterations suggests a promising adaptive strategy: starting with a smaller index set in the early stages and gradually increasing the number of snapshots as the iteration proceeds may achieve an effective balance between accuracy and computational efficiency.

\section{Conclusions}
\label{sec:conclusions}

In this paper, we explore the development of reduced order modeling (ROM) techniques for global-in-time iterative decoupled algorithms applied to Biot's consolidation model. These algorithms require repeatedly solving generalized Stokes problems at each iteration. To enhance computational efficiency, we propose a hybrid approach in which a selected subset of problems is solved using the full order model (FOM), while the remaining are approximated using ROMs constructed via proper orthogonal decomposition (POD). The selection of the FOM subset plays a critical role, as it can significantly reduce the number of required snapshots without compromising solution quality. Numerical experiments demonstrate that the proposed strategy achieves a favorable balance between efficiency and accuracy. 
The similarity observed in early-iteration errors suggests an adaptive strategy: begin with a sparse index set and gradually enrich the snapshots as iterations progress.
In future work, we aim to incorporate a posteriori error estimation techniques, as used in certified reduced basis methods, to guide the adaptive selection of the FOM subset. 
Furthermore, we will focus on extending the current framework to nonlinear poroelastic models and exploring its application to more complex, real-world poroelastic systems.

\bibliographystyle{siamplain}
\bibliography{references}

\begin{thebibliography}{10}

\bibitem{ahmed2020adaptive}
{\sc E.~Ahmed, J.~M. Nordbotten, and F.~A. Radu}, {\em Adaptive asynchronous time-stepping, stopping criteria, and a posteriori error estimates for fixed-stress iterative schemes for coupled poromechanics problems}, Journal of Computational and Applied Mathematics, 364 (2020), p.~112312.

\bibitem{akrivis2021energy}
{\sc G.~Akrivis, M.~Chen, F.~Yu, and Z.~Zhou}, {\em The energy technique for the six-step {BDF} method}, SIAM Journal on Numerical Analysis, 59 (2021), pp.~2449--2472.

\bibitem{almani2016convergence}
{\sc T.~Almani, K.~Kumar, A.~H. Dogru, G.~Singh, and M.~F. Wheeler}, {\em Convergence analysis of multirate fixed-stress split iterative schemes for coupling flow with geomechanics}, Computer Methods in Applied Mechanics and Engineering, 311 (2016), pp.~180--207.

\bibitem{almani2023convergence}
{\sc T.~Almani, K.~Kumar, and M.~F. Wheeler}, {\em Convergence analysis of single rate and multirate fixed stress split iterative coupling schemes in heterogeneous poroelastic media}, Numerical Methods for Partial Differential Equations, 39 (2023), pp.~3170--3194.

\bibitem{altmann2021semi}
{\sc R.~Altmann, R.~Maier, and B.~Unger}, {\em Semi-explicit discretization schemes for weakly coupled elliptic-parabolic problems}, Mathematics of Computation, 90 (2021), pp.~1089--1118.

\bibitem{altmann2024higher}
{\sc R.~Altmann, A.~Mujahid, and B.~Unger}, {\em Higher-order iterative decoupling for poroelasticity}, Advances in Computational Mathematics, 50 (2024), p.~111.

\bibitem{ballarin2024projection}
{\sc F.~Ballarin, S.~Lee, and S.-Y. Yi}, {\em Projection-based reduced order modeling of an iterative scheme for linear thermo-poroelasticity}, Results in Applied Mathematics, 21 (2024), p.~100430.

\bibitem{baratta2023dolfinx}
{\sc I.~A. Baratta, J.~P. Dean, J.~S. Dokken, M.~Habera, J.~HALE, C.~N. Richardson, M.~E. Rognes, M.~W. Scroggs, N.~Sime, and G.~N. Wells}, {\em {DOLFIN}x: the next generation {FE}ni{CS} problem solving environment},  (2023).

\bibitem{berkooz1993proper}
{\sc G.~Berkooz, P.~Holmes, and J.~L. Lumley}, {\em The proper orthogonal decomposition in the analysis of turbulent flows}, Annual review of fluid mechanics, 25 (1993), pp.~539--575.

\bibitem{biot1941general}
{\sc M.~A. Biot}, {\em General theory of three-dimensional consolidation}, Journal of Applied Physics, 12 (1941), pp.~155--164.

\bibitem{borregales2019partially}
{\sc M.~Borregales, K.~Kumar, F.~A. Radu, C.~Rodrigo, and F.~J. Gaspar}, {\em A partially parallel-in-time fixed-stress splitting method for {B}iot’s consolidation model}, Computers \& Mathematics with Applications, 77 (2019), pp.~1466--1478.

\bibitem{brenner1993nonconforming}
{\sc S.~Brenner}, {\em A nonconforming mixed multigrid method for the pure displacement problem in planar linear elasticity}, SIAM Journal on Numerical Analysis, 30 (1993), pp.~116--135.

\bibitem{cai2023combination}
{\sc M.~Cai, H.~Gu, P.~Hong, and J.~Li}, {\em A combination of physics-informed neural networks with the fixed-stress splitting iteration for solving {B}iot's model}, Frontiers in Applied Mathematics and Statistics, 9 (2023), p.~1206500.

\bibitem{cai2023some}
{\sc M.~Cai, H.~Gu, J.~Li, and M.~Mu}, {\em Some optimally convergent algorithms for decoupling the computation of {B}iot’s model}, Journal of Scientific Computing, 97 (2023), p.~48.

\bibitem{chaabane2018splitting}
{\sc N.~Chaabane and B.~Rivi{\`e}re}, {\em A splitting-based finite element method for the {B}iot poroelasticity system}, Computers \& Mathematics with Applications, 75 (2018), pp.~2328--2337.

\bibitem{chen2013efficient}
{\sc W.~Chen, M.~Gunzburger, D.~Sun, and X.~Wang}, {\em Efficient and long-time accurate second-order methods for the {S}tokes-{D}arcy system}, SIAM Journal on Numerical Analysis, 51 (2013), pp.~2563--2584.

\bibitem{dana2018multiscale}
{\sc S.~Dana, B.~Ganis, and M.~F. Wheeler}, {\em A multiscale fixed stress split iterative scheme for coupled flow and poromechanics in deep subsurface reservoirs}, Journal of Computational Physics, 352 (2018), pp.~1--22.

\bibitem{feng2018analysis}
{\sc X.~Feng, Z.~Ge, and Y.~Li}, {\em Analysis of a multiphysics finite element method for a poroelasticity model}, IMA Journal of Numerical Analysis, 38 (2018), pp.~330--359.

\bibitem{fischer2024more}
{\sc H.~Fischer, J.~Roth, T.~Wick, L.~Chamoin, and A.~Fau}, {\em {MOR}e {DWR}: space-time goal-oriented error control for incremental {POD}-based {ROM} for time-averaged goal functionals}, Journal of Computational Physics, 504 (2024), p.~112863.

\bibitem{fu2019high}
{\sc G.~Fu}, {\em A high-order {HDG} method for the {B}iot’s consolidation model}, Computers \& Mathematics with Applications, 77 (2019), pp.~237--252.

\bibitem{girault1979finite}
{\sc V.~Girault and P.-A. Raviart}, {\em Finite element approximation of the {N}avier-{S}tokes equations}, vol.~749, Springer Berlin, 1979.

\bibitem{gu2023iterative}
{\sc H.~Gu, M.~Cai, and J.~Li}, {\em An iterative decoupled algorithm with unconditional stability for {B}iot model}, Mathematics of Computation, 92 (2023), pp.~1087--1108.

\bibitem{gu2024crank}
{\sc H.~Gu, M.~Cai, and J.~Li}, {\em Crank-{N}icolson-type iterative decoupled algorithms for {B}iot's consolidation model using total pressure}, arXiv preprint arXiv:2409.18391,  (2024).

\bibitem{hairer2010solving}
{\sc E.~Hairer and G.~Wanner}, {\em Solving Ordinary Differential Equations \text{II}: Stiff and Differential-Algebraic Problems}, Springer Series in Computational Mathematics, Springer Berlin Heidelberg, 2010.

\bibitem{hesthaven2016certified}
{\sc J.~S. Hesthaven, G.~Rozza, and B.~Stamm}, {\em Certified reduced basis methods for parametrized partial differential equations}, vol.~590, Springer, 2016.

\bibitem{ju2020parameter}
{\sc G.~Ju, M.~Cai, J.~Li, and J.~Tian}, {\em Parameter-robust multiphysics algorithms for {B}iot model with application in brain edema simulation}, Mathematics and Computers in Simulation, 177 (2020), pp.~385--403.

\bibitem{kim2011stability}
{\sc J.~Kim, H.~A. Tchelepi, and R.~Juanes}, {\em Stability and convergence of sequential methods for coupled flow and geomechanics: Fixed-stress and fixed-strain splits}, Computer Methods in Applied Mechanics and Engineering, 200 (2011), pp.~1591--1606.

\bibitem{lee2017parameter}
{\sc J.~J. Lee, K.-A. Mardal, and R.~Winther}, {\em Parameter-robust discretization and preconditioning of {B}iot's consolidation model}, SIAM Journal on Scientific Computing, 39 (2017), pp.~A1--A24.

\bibitem{mikelic2014numerical}
{\sc A.~Mikeli{\'c}, B.~Wang, and M.~F. Wheeler}, {\em Numerical convergence study of iterative coupling for coupled flow and geomechanics}, Computational Geosciences, 18 (2014), pp.~325--341.

\bibitem{mikelic2013convergence}
{\sc A.~Mikeli{\'c} and M.~F. Wheeler}, {\em Convergence of iterative coupling for coupled flow and geomechanics}, Computational Geosciences, 17 (2013), pp.~455--461.

\bibitem{naumovich2006finite}
{\sc A.~Naumovich}, {\em On finite volume discretization of the three-dimensional {B}iot poroelasticity system in multilayer domains}, Computational Methods in Applied Mathematics, 6 (2006), pp.~306--325.

\bibitem{nevanlinna1981multiplier}
{\sc O.~Nevanlinna and F.~Odeh}, {\em Multiplier techniques for linear multistep methods}, Numerical Functional Analysis and Optimization, 3 (1981), pp.~377--423.

\bibitem{oyarzua2016locking}
{\sc R.~Oyarz{\'u}a and R.~Ruiz-Baier}, {\em Locking-free finite element methods for poroelasticity}, SIAM Journal on Numerical Analysis, 54 (2016), pp.~2951--2973.

\bibitem{phillips2007coupling}
{\sc P.~J. Phillips and M.~F. Wheeler}, {\em A coupling of mixed and continuous {G}alerkin finite element methods for poroelasticity {I}: the continuous in time case}, Computational Geosciences, 11 (2007), pp.~131--144.

\bibitem{RozzaBallarinScandurraPichi2024}
{\sc G.~Rozza, F.~Ballarin, L.~Scandurra, and F.~Pichi}, {\em Real Time Reduced Order Computational Mechanics}, SISSA Springer Series, Springer Cham, 2024.

\bibitem{siade2010snapshot}
{\sc A.~J. Siade, M.~Putti, and W.~W.-G. Yeh}, {\em Snapshot selection for groundwater model reduction using proper orthogonal decomposition}, Water Resources Research, 46 (2010).

\bibitem{suli2003introduction}
{\sc E.~S{\"u}li and D.~F. Mayers}, {\em An introduction to numerical analysis}, Cambridge university press, 2003.

\bibitem{yi2024physics}
{\sc S.-Y. Yi and S.~Lee}, {\em Physics-preserving enriched {G}alerkin method for a fully-coupled thermo-poroelasticity model}, Numerische Mathematik, 156 (2024), pp.~949--978.

\end{thebibliography}
\end{document}